\documentclass[11 pt,final]{article}
\usepackage{amsfonts,latexsym} 

\usepackage{graphicx}
\usepackage{amssymb}
\usepackage{epstopdf}
\usepackage{color}
\usepackage{aeguill}

\usepackage[color]{showkeys}

\newcounter{cst}
\def \ctel#1{C_{\refstepcounter{cst}\label{#1}\thecst}}
\def \cter#1{C_{\ref{#1}}}

\newcounter{cexp}
\def \terml#1{T_{\refstepcounter{cexp}\label{#1}\thecexp}}
\def \termr#1{T_{\ref{#1}}}

\newtheorem{defi}{Definition}[section]
\newtheorem{lemma}{Lemma}[section]
\newtheorem{remark}{Remark}[section]
\newtheorem{theo}{Theorem}[section]

\newtheorem{cor}{Corollary}[section]
\newenvironment{proof}{\noindent {\sc Proof.} }{$\square$ }

\def\dsp{\displaystyle}

\def\be{\begin{equation}}
\def\ee{\end{equation}}

\def\ba{\begin{array}{lllll}}
\def\ea{\end{array}}

\def\beqsys {\be\ba \left \{ \begin{array}{l}}
\def\eeqsys {\end{array} \right . \ea\ee }

\def\beqsysno {\be\ba \left \{ \begin{array}{l}}
\def\eeqsysno {\end{array} \right . \ea\ee}

\def\bu{\bar u}
\def\bfv{{\bf q}}
\def\centers{{\cal P}}

\def\card{{{\rm card}}}

\def\cv{K}
\def\cvv{L}

\def\drm{\rm d}

\def\disc{{\cal D}}
\def\discn{{\disc_{n}}}

\def\diam{\hbox{\rm diam}}
\def\div{{\rm div}}

\def\dkl{d_{\kl}}

\def\dr{\partial}
\def\dx{{\rm d}x}
\def\dfrontiere{{\rm d}\gamma}

\def\edge{\sigma}
\def\edges{{\cal E}}
\def\edgesint{{\cal E}_{{\rm int}}}
\def\edgesext{{\cal E}_{{\rm ext}}}
\def\edgecvcvv{K|L}
\def\edgescv{{\cal E}_K}
\def\edgescvext{{\cal E}_{K,{\rm ext}}}

\def\eps{\varepsilon}

\def\grad{\nabla}

\def\half{{\frac 1 2}}

\def\kl{{K|L}}

\def\lap{\Delta}

\def\matrices{{\cal M}_d(\R)}
\def\mcv{\meas(K)}
\def\medge{\meas(\edge)}
\def\meas{{\rm m}}
\def\mesh{{\cal M}}
\def\mkl{\meas(K|L)}

\def\n{\mathbf{n}}

\def\N{\mathbb{N}}
\def\NN{{\cal N}}

\def\normedeu(#1){\|#1\|_{L^2(\Omega)}}

\def\O{\Omega}

\def\phi{\varphi}

\def\R{\mathbb{R}}
\def\refe#1{(\ref{#1})}

\def\tends{\to}
\def\tkl{\tau_{\kl}}

\def\transedge{\tau_\edge}

\begin{document}


\title{A cell-centered finite volume approximation \\ for   second order 
partial derivative operators with full matrix \\
on  unstructured meshes in any space dimension}

\author{R. Eymard\footnote{Universit\'e de Marne-la-Vall\'ee, France,
eymard@math.univ-mlv.fr}, T. Gallou\"et\footnote{Universit\'e de Provence, France,
gallouet@cmi.univ-mrs.fr}  and 
R. Herbin\footnote{Universit\'e de Provence, France,
herbin@cmi.univ-mrs.fr}}
\maketitle

{\bf Abstract.} Finite volume methods for problems involving second
order
operators with full diffusion matrix can be used thanks to the definition
 of a discrete gradient for
piecewise constant functions on unstructured meshes satisfying an
orthogonality
condition. This discrete gradient is shown to satisfy a strong
convergence property on the interpolation of regular functions, and a weak
one on functions bounded for a discrete $H^1$ norm. To highlight
 the importance of both
properties, the convergence of the finite volume scheme on 
a homogeneous Dirichlet problem with full diffusion matrix
is proven, and an error estimate is provided. Numerical
tests show the actual accuracy of the method.

\medskip


\section{Introduction}

The approximation of convection diffusion problems in 
anisotropic media is an important issue in several engineering fields.
  Let us briefly review 
four particular  situations where
  the discretization of a nondiagonal second order operator is required:
\begin{enumerate}
\item In the case of a contaminant transported by a one-phase flow, one must
account for the diffusion-dispersion operator $\div (\Lambda\grad u)$, where
the matrix $\Lambda(x) = \lambda(x) {\rm I}_d + \mu(x) {\bfv(x)\cdot \bfv(x)^t}$ depends
on the space variable $x$ and
$\bfv(x)$ is the velocity of the fluid flow in the porous medium. The real parameter
$\lambda(x)$ corresponds
to a resulting isotropic diffusion term, including dispersion in the  directions orthogonal to
the flow, and the real parameter $\mu(x)$ to an additional
diffusion in the direction of the flow \cite{letem}. The term $\bfv(x)$ is then given
by  $\bfv(x) = K(x) \grad p(x)$, where  $p(x)$ is a pressure and
$K(x)$ another nondiagonal matrix (the absolute permeability matrix, 
depending on the geological layers), and
satisfies the incompressibility equation $\div \bfv(x) = 0$. In this coupled
problem, one must simultaneously compute this pressure and the contaminant concentration $u(x)$.

\item In the study
of undersaturated flows in porous media (for example, air-water flows), two equations
of conservation have to be solved, associated with two unknowns,   pressure
and   saturation.
These equations include nonlinear hyperbolic and degenerate parabolic terms with
respect to the saturation unknown. As in the preceding case, one must
discretize such terms as  $\div \bfv(x) = \div (K(x) \grad p(x))$, where 
again $K(x)$ is a nondiagonal matrix
depending on the geological layers.

\item In the case of the compressible Navier-Stokes equations,  
one has to discretize the viscous forces operator, which can be written  under the form
$a \lap {\bf u} + b \grad \div {\bf u}$
($a$ and $b$ are deduced from the dynamic viscosity coefficients and ${\bf u}$ is
the fluid velocity). In this problem, the term $\grad \div {\bf u}$ involves
all the cross derivatives $\dr_{ij}^2 {\bf u}$.

\item
Some problems arising in financial mathematics lead to anisotropic
diffusion equations in high-dimensional domains (dimension equal to
5 or more for example). 
Under some assumptions on financial markets \cite{LL}, 
the price of a European or an American option is obtained
by solving a linear or nonlinear partial differential equation, 
 involving the second order anisotropic diffusion matrix $\Lambda = \Sigma \Sigma^t$, where 
$\Sigma$ is a real matrix.

\end{enumerate}

All these cases involve a term under the form
$\div (\Lambda\grad u)$, where $\Lambda$ is 
a (generally) nondiagonal matrix depending on the space variable and $u$ 
is a function of 
the space variable in steady problems, and of the space and time variables in
transient problems. 
Finite element schemes are known to allow for an 
easy discretization of such
a term on triangular or tetrahedral meshes \cite{putti}.
However, in engineering situations such as the ones
described above, one  also has to discretize convection and reaction terms, and avoid 
numerical instabilities. Unfortunately,  finite element methods (and more generally centered schemes)
 are known to generate instabilities 
on coarse grids, although some cures may be 
proposed, see \cite{f,angermann}; therefore
a great many numerical codes \cite{aav1,aav2,f,pt,jt} 
use finite volume or finite volume - finite element type schemes, 
which allow the implementation of discretization techniques (such as the classical upwind
schemes)
which prevent  the apparition of instabilities.
Let us also note that  finite volume schemes are known for their simplicity of
implementation, particularly so when discretizing coupled systems of 
equations of various nature. 

\medskip 

Besides, a thorough mathematical analysis has now been improved,
showing that finite volume methods are well suited and convergent 
for a simple convection diffusion equation  in the case where 
$\Lambda(x) = \lambda(x) \ {\rm I}_d$.
Indeed, this analysis has been completed  (see
\cite{vf4},  \cite{mishev}, \cite{ghv}, \cite{book}) in the case 
of grids (called admissible in the sense of \cite{book}, see also Definition \ref{adisc} below) 
satisfying an orthogonality condition:
the line joining two cell centers is orthogonal to the
interface between the two cells,  thus  ensuring
a consistency property when approximating 
the normal flux at the cell interface by centered finite differences. 
Some examples of such admissible grids are
the Delaunay triangular meshes or tetrahedral meshes,
  rectangular or parallelepipedic meshes in 2 or 3 dimensions, and 
the Vorono\"{\i} meshes in any dimension.

\medskip

But the situation is quite different in the case where  
the condition $\Lambda(x) = \lambda(x) \ {\rm I}_d$ 
  no longer holds:
only few of the actual discretization methods used for handling nondiagonal second
order terms on finite volume grids meet a full mathematical
analysis of stability or convergence. 
Let us briefly review some of them. A first one, in the case where 
$\Lambda(x) = \lambda(x) \ M$, where $M$ is a symmetric  positive definite 
matrix,  consists in
adapting the above orthogonality condition by stating that the line joining two cell centers 
is orthogonal to the
interface between the two cells with respect to the dot product induced by the
matrix $\Lambda^{-1}$. Indeed, it is also possible to consider the case 
where $M$ depends on the discretization cell, by using, in each cell, 
 the orthogonal bisectors
for the metric induced by $M^{-1}$ (see \cite{fvca1} and \cite{book} 
section 11 page 815).
In the case of triangular grids, this yields a 
well defined scheme under some restriction on the allowed anisotropy 
for a given geometry, 
since  the cell center is chosen
as the  intersection of the orthogonal bisectors of the triangle for 
the metric defined by  $M^{-1}$.   Another method consists 
in defining the finite volume method as a dual method
to a finite element one (for example, a P1 finite element \cite{letem} 
or a Crouzeix-Raviart one,
see e.g. \cite{voralhik}).

\medskip

Another possibility to derive a  finite volume scheme on problems including
anisotropic diffusion is to construct  a 
local discrete gradient, allowing to get, at each edge $\edge$ of the mesh,
a consistent approximate value for the flux
$\int_\edge (\Lambda(x) \grad u(x))\cdot \n_{\edge} \dfrontiere(x)$
involved in the finite volume scheme ($\n_{\edge}$ is a unit
vector normal to the edge $\edge$, and $\dfrontiere(x)$ is the
$d-1$ Lebesgue measure on the edge $\edge$ ).
 In two space dimensions, such a scheme was introduced 
in  \cite{coudiere} 
on arbitrary meshes,
but the proof of convergence was only possible on meshes close to
parallelograms. Still in 2D,  a technique using dual meshes   is 
introduced in \cite{hermeline,omnes}, which generalizes the idea of
\cite{nic1,nic2} for div-curl problems to meshes with no orthogonality conditions; however 
the use of a dual mesh renders the scheme computationally expensive; moreover
it does not seem  to be easily extended to 3D.
In \cite{convgrad}, we used
Raviart-Thomas shape functions, generalized to the case of any admissible mesh 
(again in the sense precised
 of \cite{book}, see also Definition \ref{adisc} below), in order to
define a discrete gradient for piecewise constant functions. 
The strong convergence of this discrete gradient was then
shown in the case of the  elliptic equation $-\lap u = f$. A drawback
of this definition was the difficulty to find an approximation of these generalized shape functions
in other cases than triangles or rectangles. 

\medskip

We therefore propose in this paper a new cheap and simple 
method of constructing a 
discrete gradient for
a piecewise constant function, 
on arbitrary admissible meshes in any space dimension (this method has been 
first introduced in \cite{crasvf100}). We prove that the discrete gradients of 
any sequence of 
piecewise constant functions converging to some $u \in H^1_0(\Omega)$ weakly converges
to $\grad u $ in $L^2(\Omega)$.  Moreover, the discrete gradient is shown to 
be consistent, in the sense that it  satisfies 
 a strong convergence property on the
interpolation
 of regular function. In order to show the efficiency of this approximation method,
we use this  discrete gradient to design a scheme for the approximation of
the weak solution $\bar u$ of the following diffusion  problem with 
full anisotropic tensor:
\be\ba\dsp
 - \div ( \Lambda \grad \bar u)  = f \hbox{ in } \Omega, \\
 \bar u=0 \hbox{ on } \partial \Omega,
\ea\label{ellgen}\ee
under the following assumptions:
\be
\O \mbox{ is an open bounded connected polygonal subset of }\R^d,
\ d\in\N^\star,
\label{hypomega}\ee

\be
\ba
\Lambda \hbox{ is a measurable function from }  \Omega \hbox{ to } \matrices,
\\  
\hbox{ where  }\matrices 
\mbox{ denotes the set of }  d\times
d  \mbox{  matrices, }
\\
\mbox{ such that for a.e. } x \in \Omega, \Lambda(x) 
\mbox{ is symmetric,}
\\ \mbox{ and 
the set of its eigenvalues is included in } [\alpha(x),\beta(x)]
\\
\mbox{ where } \alpha,\beta\in L^\infty(\O) \mbox{ are  such that }
\\
0<\alpha_0\le\alpha(x)\le\beta(x) \mbox{ for a.e. } x\in\O,
\ea
\label{hyplambda}\ee
and 
\be
f \in L^2(\O).
\label{hypfg}\ee
We give the classical weak formulation in the following definition.
\begin{defi}[Weak solution] \label{weaksol} 
Under hypotheses \refe{hypomega}-\refe{hypfg}, we say that $\bar u$ is a weak solution 
of \refe{ellgen} if
\be\left\{\ba
\bar u \in H^1_0(\O),\\
\dsp  \int_\O \Lambda(x)\grad \bar u(x)\cdot\grad v(x) 
\dx = \int_\O f(x)  v(x) \dx,
 \ \ \ \forall v \in H^1_0(\O).
\ea\right.\label{ellgenf}\ee
\end{defi}
\begin{remark}
For the sake of clarity,  we restrict ourselves here to the numerical analysis of Problem \refe{ellgen}, 
however, the present   analysis readily extends to convection-diffusion-reaction problems
and coupled problems. Indeed, we emphasize that   proofs of convergence 
or error estimate
can easily be adapted to such situations, since the discretization methods
of all these terms are independent of one another, and the treatment of convection and 
reaction term is well-known exact(see  \cite{ghv} or \cite{book}).
\end{remark}
\medskip

The outline of this paper is the following. In Section  \ref{secfvslin}, we present
the method for approximating the gradient of a piecewise constant function, and we
show some functional properties which help to understand why the present definition
of a gradient is well suited for second order diffusion problems. In Section \ref{seccvglin},
we present the finite volume scheme for Problem \refe{ellgen}, and we show the strong convergence of the discrete solution and of its discrete gradient.
In Section \ref{errest}, we give an error estimate for Problem  \refe{ellgen}, and 
we illustrate this study by some numerical examples in  Section \ref{numex}.
Some short conclusions are drawn in Section \ref{conclu}.

\section{A discrete gradient for piecewise constant functions}\label{secfvslin}

We present in this section a method for the approximation of 
the gradient of  piecewise constant functions, in the case of grids satisfying
some orthogonality condition as defined below.

\subsection{Admissible discretization of $\O$}

\begin{figure}[ht]
\begin{center}
\input{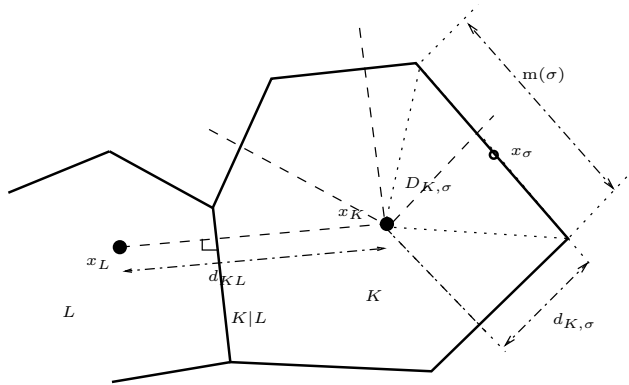}
\end{center}
\caption{Notations for a control volume $K$ in the case $d=2$}\label{fig_maille}
\end{figure}

We first present the following notion of admissible discretization, which
is taken in \cite{book}. The notations are summarized in Figure \ref{fig_maille}
for the particular case $d=2$ (we recall that the case $d\ge 3$ is considered as well).

In the following definition, we shall say that a bounded subset of $\R^d$
is polygonal if its boundary is included in the union of a finite number of 
hyperplanes. 
\begin{defi}\label{adisc}{\bf [Admissible discretization]}
Let $\O$ be an open bounded  polygonal 
subset of $\R^d$, and 
$\dr \O =  \overline{\O}\setminus\O$ its boundary.
An admissible finite volume discretization of $\O$,
denoted by $\disc$, is given
by $\disc=(\mesh,\edges,\centers)$, where:

\begin{itemize}
\item $\mesh$ is a finite family of
non empty open polygonal convex disjoint subsets of
$\O$ (the ``control volumes'') such that
$\overline{\O}= \dsp{\cup_{K \in \mesh} \overline{K}}$.
For any $K\in\mesh$, let
$\dr K  = \overline{K}\setminus K$ be the boundary of $K$
and $\mcv>0$ denote the measure of $K$.

\item $\edges$ is a finite family of disjoint subsets of
$\overline{\O}$ (the ``edges'' of the mesh), such that,
for all $\edge\in\edges$, there exists
a hyperplane $E$ of $\R^d$ and $K\in\mesh$
with $\overline{\edge} = \dr K \cap E$ and
$\edge$ is a non empty open subset of $E$.
We then denote by $m_\edge>0$ the $(d-1)$-dimensional measure of
$\edge$. We assume that,
for  all $K \in \mesh$, there exists  a subset $\edgescv$ of $\edges$
such that $\dr K  =
\dsp{\cup_{\edge \in \edgescv}}\overline{\edge} $.
It then results from the previous hypotheses that,
for all $\edge\in\edges$,
either $\edge\subset \dr\O$ or there exists
 $(K,L)\in \mesh^2$ with $K \neq L$ such that
$\overline{K} \cap \overline{L} = \overline{\edge}$;
we denote in the latter case $\edge = K|L$.

\item $\centers$ is a family of points of $\O$
indexed by $\mesh$, denoted by $\centers = (x_K)_{K \in \mesh}$.
The coordinates of $x_K$ are denoted by $x^{(i)}_K$, $ i = 1,\ldots,d.$
The family $\centers$ is  such that, for all $K \in \mesh$,
$x_K \in K$. Furthermore, for all $\edge\in\edges$
such that there exists  $(K,L)\in \mesh^2$ with $\edge = K|L$,
it is assumed  that  the straight line $(x_K,x_L)$
going through $x_K$ and $x_L$ is orthogonal to $K|L$.
For all $K \in \mesh$ and all $\edge \in \edgescv$,
let $z_\edge$ be the orthogonal projection of $x_K$ on $\edge$.
We suppose that $z_\edge\in\edge$ if $\edge \subset \partial \Omega$.

\end{itemize}
 
\end{defi}

The following notations are used.
The size of the discretization is defined by:

\[
h_\disc = \sup\{\hbox{\rm diam}(K), K\in \mesh\}.
\]

For all $K \in \mesh$ and $\edge \in \edgescv$, we denote
by $\n_{K,\edge}$ the unit vector normal to $\edge$ outward to $K$.
We denote by $d_{K,\edge}$
the Euclidean distance between $x_K$ and $\edge$. We then define
\[
\tau_{K,\edge} = \frac {m_\edge} {d_{K,\edge}}.
\]

The set of interior (resp. boundary) edges is denoted by $\edgesint$
(resp. $\edgesext$), that is $\edgesint = \{\edge \in \edges;$ $\edge \not \subset
\partial \O \}$ (resp.  $\edgesext = \{\edge \in \edges;$ $\edge \subset
\partial \O \}$). For all
$K\in\mesh$, we denote by $\NN_K$ the subset of $\mesh$ of the 
neighbouring control volumes, and we denote by
$\edgescvext = \edgescv\cap\edgesext$.
For all  $\edge\in\edgesint$, let $K,L\in\mesh$ be such that $\edge = K|L$;
 we define by $d_{K \vert L}$ the Euclidean distance between $x_K$ and $x_L$, by
$\n_{KL}$ the unit normal vector to $K|L$ from $K$ to $L$, and we set
\be
\tau_{\edge} = \frac {m_\edge} {d_{K \vert L}}.
\label{tauedgeint}\ee
For all $\edge\in\edgesext$, let $K\in\mesh$ be such that $\edge\in\edgescv$; we define
\be 
\tau_{\edge} = \tau_{K,\edge}.
\label{tauedgeext}\ee
For all $K\in\mesh$ and $\edge\in\edgescv$, we define
\[
D_{K,\edge}=\{ t x_K +(1-t) y, t\in (0,1), \ y\in \edge\},
\]
For all $\edge\in\edgesint$, let $K,L\in\mesh$ be such that $\edge = K|L$;
we set $D_\edge = D_{K,\edge}\cup D_{L,\edge}$.
For all $\edge\in\edgesext$, let $K\in\mesh$ be such that $\edge\in\edgescv$;
we define $D_\edge = D_{K,\edge}$.

For all $\edge\in\edges$, we define
\be
 x_{\edge} = \frac {1} {\medge} \int_{\edge} x \ {\rm d}\gamma(x).
\label{xedge}\ee 
We shall measure the regularity of the mesh through the 
 function $\theta_\disc$ defined by
\be
\theta_\disc  = \inf\left\{\frac {d_{K,\edge}} {\hbox{\rm diam}(K)},
 K\in \mesh,\ \edge\in\edgescv \right\}.
\label{reguld}
\ee
\begin{defi}\label{defHdisc}
Let $\Omega$ be an open bounded polygonal subset of $\R^d$, 
and $\disc $ an admissible discretization of $\Omega$ in the sense 
of Definition \refe{adisc}. We define 
 $H_\disc$  as the set of functions $u\in L^2(\O)$ which are
constant in each control volume. For $u\in H_\disc$, we denote by
$u_\cv$ the constant value of $u$ in $\cv$.
We define the interpolation operator
$P_\disc~:~C(\overline \O)\to H_\disc$, by $\bu\mapsto P_\disc \bu$
such that
\be
P_\disc \bu (x) = \bu(x_\cv) \hbox{ for a.e. } x \in \cv,\ \forall \cv \in \mesh.
\label{Pdisc}
\ee
For $(u,v)\in (H_\disc)^2$ and for any function $\alpha\in L^\infty(\O)$, 
we introduce the following symmetric bilinear form:
\be
[u,v]_{\disc,\alpha} = \sum_{K|L\in\edgesint} \tkl \alpha_{K|L} (u_L-u_K)(v_L - v_K) +
\sum_{K \in \mesh} \sum_{\edge\in \edgescvext} \tau_{K,\edge} \alpha_{\edge} u_K v_K, 
\label{amudis}\ee 
where we set
\be 
\alpha_{\edge} = \frac {1} {\meas(D_\edge)} \int_{D_\edge} \alpha(x)\dx,\ \forall \edge\in\edges.
\label{alphaedge}\ee 
\end{defi}
 
\begin{remark}
One could also take, for $\alpha_{\edge}$, the 
harmonic averaging of the values in $K$ and $L$ when $\sigma = K|L$.
\end{remark}
We then define a norm in $H_\disc$  (thanks to the discrete Poincar\'e
inequality \refe{poindis} given below) by
\[
\Vert u\Vert_{\disc} = \left( [u,u]_{\disc,1} \right)^{1/2} 
\]
(where 1 denotes the constant function equal to 1).
Indeed, the discrete Poincar\'e inequality writes (see \cite{book}):
\be
\Vert w \Vert_{L^2(\O)} \le \diam(\O) \Vert w\Vert_{\disc}, \ \forall w\in H_\disc.
\label{poindis}\ee
Let us now give a relative compactness result, which is also partly stated in some
other papers concerning finite volume methods \cite{book}, \cite{convpardeg}.
\begin{lemma}[Relative compactness in $L^2(\O)$]\label{discpct}
Let $\O$ be an open bounded connected polygonal subset of $\R^d$,
$d\in\N^\star$ and let $(\disc_n,u_n)_{n\in\N}$ be a sequence
such that, for all $n\in\N$, $\disc_n$ is an
admissible finite volume discretization
 of $\O$ in the sense of Definition
\ref{adisc} and  $u_n \in H_{\disc_n}(\O)$
(cf Definition \ref{defHdisc}). Let us assume that
$\lim_{n\tends\infty} h_\discn  = 0$, and
that there exists $\ctel{u1}>0$ such that 
 $\Vert u_n\Vert_{\disc_n} \le \cter{u1}$, for all  
 $n\in\N$.

Then there exists a subsequence of $(\disc_n,u_n)_{n\in\N}$,
again denoted $(\disc_n,u_n)_{n\in\N}$, and $\bu\in H^1_0(\O)$
such that $u_n$ tends to $\bu$ in $L^2(\O)$ as $n\tends+\infty$, and the inequality
\be
\int_\O |\grad\bu(x)|^2 \dx \le \liminf_{n\tends\infty} \Vert u_n\Vert_{\disc_n}^2
\label{liminf}
\ee
holds. Moreover, for all function $\alpha\in L^\infty(\O)$, we have
\be
\lim_{n\tends\infty} [u_n,P_{\disc_n}\varphi]_{\disc_n,\alpha} = 
\int_\O \alpha(x) \grad\bu(x)\cdot\grad\varphi(x)\dx,\ \forall \varphi\in {\rm C}^\infty_c(\O).
\label{consdis}\ee
\end{lemma}
\begin{proof}
The proof of the existence of the subsequence
again denoted $(\disc_n,u_n)_{n\in\N}$, and of 
$\bu\in H^1_0(\O)$ such that $u_n$ tends to $\bu$
in $L^2(\O)$ as $n\tends\infty$, is given in \cite{book}.  Assertion \refe{liminf}
was proven in \cite{convpardeg} (Lemma 5.2).
Let us first show \refe{consdis} in the case $\alpha\in {\rm C}^1(\bar\O)$. Let $\varphi\in {\rm C}^\infty_c(\O)$. 
Defining, for all $n\in\N$,
$\terml{xa1}^{(n)} = - \int_\O u_n(x) 
\div(\alpha(x) \grad\varphi(x))\dx$,
we get that 
\[
\lim_{n\to\infty} \termr{xa1}^{(n)} = - \int_\O \bu(x) 
\div(\alpha(x) \grad\varphi(x))\dx = \int_\O \alpha(x) \grad\bu(x)\cdot\grad\varphi(x)\dx.
\]
We consider a value $n$ sufficiently large such that for all $K\in\mesh_n$ and $x\in K$,
if $\varphi(x)\neq 0$ then $\dr K\cap \dr\O = \emptyset$. 
Defining $\terml{xa2}^{(n)} = [u_n,P_{\disc_n}\varphi]_{\disc_n,\alpha} - \termr{xa1}^{(n)}$, we obtain
\[
\termr{xa2}^{(n)} = \sum_{\sigma\in\edgesint,\ \sigma = K|L} \mkl  (u_L-u_K) R_{KL},
\]
with
\[
R_{KL} = \alpha_{K|L}\frac{\varphi(x_L) - \varphi(x_K)} {\dkl} - \int_{K|L} \alpha(x) 
\grad\varphi(x)\cdot \n_{KL}
\dfrontiere(x),\ \forall K\in\mesh,\ \forall L\in\NN_K .
\]
Since there exists some real value $\ctel{aphi}$, which does not depend on $\disc_n$, such that
$|R_{KL}|\le \cter{aphi}h_\discn $, we conclude in a similar way as in \cite{book} that
$\lim_{n\to\infty} \termr{xa2}^{(n)} = 0$, which gives \refe{consdis} in this case.
Let us now consider the general case $\alpha\in L^\infty(\O)$.
Let $\eps>0$ be given.
We first choose a function $\tilde\alpha\in {\rm C}^1(\bar\O)$ such that 
$\Vert \alpha - \tilde\alpha\Vert_{L^2(\O)}\le\eps$.
Then we have, for all $n\in\N$, using the Cauchy-Schwarz inequality,
\[\ba
\dsp
\left([u_n,P_{\disc_n}\varphi]_{\disc_n,\tilde\alpha}-[u_n,P_{\disc_n}\varphi]_{\disc_n,\alpha}\right)^2
\le
& \dsp \sum_{K|L \in\edgesint} \tkl(\tilde\alpha_{KL}-\alpha_{KL})^2 |\varphi(x_L)-\varphi(x_K)|^2 \\

& \dsp \times \sum_{K|L \in\edgesint} \tkl |u_L - u_K|^2
\ea\]
and therefore, setting $\ctel{phiinf} = \Vert\grad\varphi\Vert_{L^\infty(\O)}$,
the properties $|\varphi(x_L)-\varphi(x_K)|\le \cter{phiinf} d_{K \vert L}$ and
$\mkl\dkl = d \ \meas(D_{K|L})$ lead to
\[\dsp
\left([u_n,P_{\disc_n}\varphi]_{\disc_n,\tilde\alpha}-[u_n,P_{\disc_n}\varphi]_{\disc_n,\alpha}\right)^2\le
d\ \cter{phiinf}^2 \Vert \alpha - \tilde\alpha\Vert_{L^2(\O)}^2 \cter{u1}\le 
d\ \cter{phiinf}^2 \eps^2 \cter{u1}.
\]
In the same manner, we get
\[
\left(\int_\O \tilde\alpha(x) \grad\bu(x)\cdot\grad\varphi(x) \dx - 
\int_\O \alpha(x) \grad\bu(x)\cdot\grad\varphi(x) \dx
\right)^2\le \cter{phiinf}^2 \eps^2 \Vert \grad\bu\Vert_{L^2(\O)^d}^2.
\]
Since $\tilde\alpha\in {\rm C}^1(\O)$, we can apply \refe{consdis}, proven above for such a function.
It then  suffices to choose $n$ large enough such that
\[
\left|[u_n,P_{\disc_n}\varphi]_{\disc_n,\tilde\alpha}-
\int_\O \tilde\alpha(x) \grad\bu(x)\cdot\grad\varphi(x) \dx\right| \le\eps,
\]
to prove that
\[
\left|[u_n,P_{\disc_n}\varphi]_{\disc_n,\alpha}-\int_\O \alpha(x) 
\grad\bu(x)\cdot\grad\varphi(x) \dx\right| \le \ctel{u2}\ \eps,
\]
where the real $\cter{u2}>0$ does not depend on $n$. This concludes the
proof of \refe{consdis} in the general case.
\end{proof}

\subsection{Definition of a discrete gradient}

We now define a discrete gradient for piecewise constant functions on 
an admissible discretization.

\begin{defi}[Discrete gradient]\label{espdis}
Let $\O$ be an open bounded connected polygonal subset of $\R^d$,
$d\in\N^\star$. Let 
$\disc=(\mesh,\edges,\centers)$ be an
 admissible finite volume discretization of $\O$ in the sense of Definition
\ref{adisc}. 
Let us define, for all $K\in\mesh$, for all $L\in\NN_K$,
\be\ba\dsp
A_{K,L} =  \tkl( x_{K|L} - x_K),
\ea\label{defaint}\ee
and for all $\edge\in\edgescvext$, we define
\be\ba\dsp
A_{K,\edge} =  \tau_\edge(x_{\edge} - x_K).
\ea\label{defaext}\ee
We define the discrete gradient  
$\grad_{\disc}~: H_\disc\to H_\disc^d$, for any $u\in H_\disc$, by:
\[\ba\dsp
\grad_{\disc}u(x)  & = (\grad_{\disc}u)_K \\&  = 
\dsp\frac 1 {\mcv}
\left(\sum_{L\in\NN_K}  A_{K,L} \ (u_L-u_K)-
\sum_{\edge\in \edgescvext}  A_{K,\edge}\ u_K\right), \\
 &\mbox{ for a.e. }x \in K,\
\forall  K\in\mesh.
\ea\]
\end{defi}

Let us first state a bound for the $L^2(\O)^d$ norm of the discrete gradient of
any element of $H_\disc$.
\begin{lemma}[Bound for $\grad_\disc u$]\label{boundnorm}
Let $\O$ be an open bounded connected polygonal subset of $\R^d$,
$d\in\N^\star$, let  $\disc$ 
be an admissible finite volume discretization of $\O$ in the sense of  Definition
\ref{adisc} and let $\theta\in (0,\theta_\disc ]$.
Then, there exists $\ctel{C}$, only depending on $d$
and $\theta$, such that, for all $u \in H_\disc$:
\be
\Vert \grad_\disc u \Vert_{L^2(\O)^d}  \le \cter{C} \Vert u \Vert_\disc.
\label{mn}
\ee
\end{lemma}

\begin{proof}
Let $u \in H_\disc$. Let us denote, for all $\cv \in \mesh$, $\cvv\in\NN_K$ and $\sigma = K|L$,
$\delta_{K,\sigma}u = u_L - u_K$, and for $\edge\in \edgescvext$,
$\delta_{K,\sigma}u =  - u_K$. Then Definition \refe{amudis} leads to 
\[
\Vert u \Vert_\disc^2 = \sum_{K \in \mesh} \left(\half \sum_{L\in\NN_K} \tkl (\delta_{K,K|L}u)^2 +
\sum_{\edge\in \edgescvext} \tau_{\edge} (\delta_{K,\edge}u)^2\right), 
\]
and Definition \refe{espdis} leads, for a given $\cv \in \mesh$, to
\[
\dsp \mcv (\grad_\disc u)_\cv  = \sum_{\edge\in\edgescv}  \transedge (x_\edge -x_\cv) \delta_{K,\sigma}u.
\]
Using the Cauchy-Scharwz inequality, we obtain
\[
\mcv^2 \vert (\grad_\disc u)_\cv \vert ^2 \le \sum_{\edge \in \edgescv} \transedge \vert x_\edge
-x_\cv \vert^2 \sum_{\edge \in \edgescv} \transedge (\delta_{K,\edge}u)^2,
\]
and, since, for $\edge \in \edgescv$,  one has $\vert x_\edge-x_\cv \vert = {\rm d}(x_\edge,x_\cv) \le
\frac {d_{\cv,\edge}} {\theta}$,
\be
\mcv^2 \vert (\grad_\disc u)_\cv \vert ^2 \le \sum_{\edge \in \edgescv} 
\frac 1 {\theta^2} \medge d_{\cv,\edge}  \sum_{\edge \in \edgescv} \transedge (\delta_{K,\sigma}u)^2.
\label{inter}
\ee
Since $\sum_{\edge \in \edgescv} \medge d_{\cv,\edge}=d\ \mcv$, \refe{inter}
gives:
\[
\mcv \vert (\grad_\disc u)_\cv \vert ^2 \le \frac d {\theta^2} 
\sum_{\edge \in \edgescv} \transedge (\delta_{K,\sigma}u)^2.
\]
Summing over $\cv \in \mesh$, we get
\[
\Vert \grad_\disc u \Vert^2_{L^2(\O)^d} \le 2  \frac d {\theta^2} \Vert u \Vert_\disc^2.
\]
which gives \refe{mn} with $\cter{C}= (\frac {2d} {\theta^2})^{\half}$. 
\end{proof}

We now state a weak convergence property for the discrete gradient.

\begin{lemma}[Weak convergence of the discrete  gradient]\label{wkcvg}~

Let $\O$ be an open bounded connected polygonal subset of $\R^d$,
$d\in\N^\star$, let $\disc$ 
be an admissible finite volume discretization of $\O$ in the sense of  Definition
\ref{adisc}. We assume that there exist $u_\disc\in H_\disc$ and
a function $\bu\in H^1_0(\O)$ such that $u_\disc$ tends to $\bu$ in $L^2(\O)$  
as $h_\disc $ tends to 0 while $\Vert u_\disc\Vert_{\disc}$ remains bounded. 
Then $\grad_{\disc} u_\disc$ weakly tends to $\grad\bu$ in $L^2(\O)^d$ as $h_\disc \tends 0$.
\end{lemma}
\begin{proof}
Let $\varphi\in {\rm C}^\infty_c(\O)$. 
We assume that $h_\disc $ is small enough to ensure that
for all  $K\in\mesh$ and $x\in K$, if
$\varphi(x) \neq 0$ then $\edges_{\cv, {\rm ext}} = \emptyset$. 
The expression $\terml{A}^\disc$, defined by
\[
\termr{A}^\disc = \int_\O P_\disc\varphi(x) \grad_{\disc} u_\disc(x) \dx,
\]
satisfies, using \refe{defaint},
\[
\termr{A}^\disc = \sum_{K|L\in\edgesint} \tkl (u_L - u_K) \left(
(x_{K|L} - x_K)\varphi(x_K) +(x_L - x_{K|L})\varphi(x_L)\right),
\]
where we denote, for the sake of simplicity, $u_K = (u_\disc)_K$ for all $K\in\mesh$.
We thus get  $\termr{A}^\disc = \terml{A1}^\disc + \terml{A2}^\disc$ with
\[
\termr{A1}^\disc = \sum_{K|L\in\edgesint} \tkl (u_L - u_K) 
(x_L - x_K)\frac {\varphi(x_K) + \varphi(x_L)} 2
\]
and 
\[
\termr{A2}^\disc = \sum_{K|L\in\edgesint} \tkl (u_L - u_K) 
(x_{K|L} - \frac {x_L + x_K} 2 ) (\varphi(x_L) -  \varphi(x_K)).
\]
Thanks to the Cauchy-Schwarz inequality, we get
\[
(\termr{A2}^\disc)^2 \le \sum_{K|L\in\edgesint} \tkl (u_L - u_K)^2
\sum_{K|L\in\edgesint} \tkl (\varphi(x_L) -  \varphi(x_K))^2
|x_{K|L} - \frac {x_L + x_K} 2 |^2 .
\]
Since $|x_{K|L} - \frac {x_L + x_K} 2 |\le \half |x_{K|L} - x_L| + \half|x_{K|L} - x_K| \le h_\disc $,
there exists $\ctel{cwkcvg}>0$, depending on $d$, $\O$ and $\varphi$ such that, 
\[
(\termr{A2}^\disc)^2 \le \Vert u_\disc\Vert_{\disc}^2 \cter{cwkcvg} h_\disc ^2 \meas(\O),
\]
and therefore we get
\[
\lim_{h_\disc \tends 0} \termr{A2}^\disc = 0. 
\]
We then compare $\termr{A1}^\disc$ with
\[
\terml{B}^\disc = -\int_\O u_\disc(x)\grad\varphi(x) \dx
=\sum_{K|L\in\edgesint} (u_L - u_K) \int_{K|L} \varphi(x)\n_{K,L} {\rm d}\gamma(x).
\]
Since 
\[
\n_{K,L} = \frac {x_L - x_K} {\dkl}
\]
and since
\[
\left| \frac 1 {\mkl} \int_{K|L} \varphi(x)  {\rm d}\gamma(x) -
\frac {\varphi(x_K) + \varphi(x_L)} 2 \right| \le \Vert\grad \varphi\Vert_{L^\infty(\O)}\  h_\disc ,
\]
we get, thanks to the Cauchy-Schwarz inequality,
\[
\lim_{h_\disc \tends 0} (\termr{A1}^\disc - \termr{B}^\disc)^2 = 0. 
\]
Since 
\[
\lim_{h_\disc \tends 0} \termr{B}^\disc = -\int_\O \bu(x) \grad\varphi(x) \dx = 
\int_\O \varphi(x) \grad\bu(x) \dx,
\]
we have thus proven, thanks to the density of ${\rm C}^\infty_c(\O)$ in $L^2(\O)$, 
the weak convergence of $\grad_{\disc} u_\disc$ to $\grad\bu(x)$ as $h_\disc \tends 0$.
This completes the proof of the lemma.
\end{proof}

We now study, for a regular function $\varphi$, the strong convergence
of the discrete gradient $\grad_\disc P_\disc\varphi$ to $\grad\varphi$.
This study uses the following lemma.
\begin{lemma}\label{agbnd} 
Let $\O$ be an open bounded connected polygonal subset of $\R^d$,
$d\in\N^\star$, let $\disc$ 
be an admissible finite volume discretization of $\O$ in the sense of  Definition
\ref{adisc}. Then we have
\be
v= \frac 1 {\mcv} \sum_{\edge\in\edgescv} \medge (x_\edge-x_0)\ 
(\n_{\cv,\edge}\cdot v)
,\ \forall \cv\in\mesh,\
\forall x_0 \in \R^d,\ \forall v \in \R^d.
\label{chxs}
\ee
\end{lemma}

\begin{proof}
For any $\cv\in\mesh$, we denote, for a.e. $x\in\dr\cv$, 
by  $\n_{\dr\cv}(x)$ the normal vector to $\partial \cv$ at the point $x$ outward $\cv$.
Let $v$ and $w \in \R^d$ be given. We have, considering vectors as $d\times 1$ matrices,
and denoting by $w^t$ the transposed $1\times d$ matrix of $w$,
\[
\ba
\dsp w^t\ \left(\int_{\partial \cv} (x-x_0) \n_{\cv}^t(x) \dfrontiere(x) 
\right)\ v=
\int_{\partial \cv} w^t\ (x-x_0)\ \n_{\cv}^t(x)\ v\ \dfrontiere(x)=
\\
\dsp \hfill \int_{\partial \cv} w^t\ (x-x_0)\ v^t \n_{\cv}(x) \dfrontiere(x)= 
\int_{\partial \cv} (v\ (x-x_0)^t\ w) \cdot \n_{\cv}(x)\dfrontiere(x)=
\\
\hfill \dsp 
\int_\cv \div (v\ (x-x_0)^t\  w) dx= \mcv\ v^t\ w.
\ea
\]
This gives  \refe{chxs}. \end{proof}

\begin{lemma}[Consistency  property of the discrete gradient]\label{cg}
Let $\O$ be an open bounded connected polygonal subset of $\R^d$,
$d\in\N^\star$,
let  $\disc$ be an admissible finite volume discretization in the sense of  Definition
\ref{adisc} and let $\theta\in (0,\theta_\disc ]$. 
Let $\bu \in {\rm C}^2(\overline \Omega)$ be such that $\bu=0$ on the boundary of $\O$.
Then, there exists $\ctel{cstcg}$, only depending on $\Omega$,
$\theta$ and $\bu$, such
that:
\be
\Vert \grad_\disc P_\disc \bu - \grad \bu \Vert_{L^2(\O)^d} \le \cter{cstcg}\ h_\disc .
\label{pdumu}
\ee
(Recall that $P_\disc$ is defined by \refe{Pdisc} and $\grad_\disc$ in Definition \ref{espdis}.)
\end{lemma}

\begin{proof} From Definition \ref{espdis} and \refe{Pdisc}, we can write
for any $\cv \in \mesh$
\be
\dsp \mcv  (\grad_\disc P_\disc\bu)_\cv  = \sum_{L\in\NN_K}  \tkl (x_{K \vert L}-x_\cv) (\bu(x_L)-\bu(x_K))-
\sum_{\edge\in \edgescvext}  \transedge (x_\edge - x_\cv) \bu(x_K) .
\label{cons}
\ee
Let $(\grad \bu)_\cv$ be the mean value of $\grad \bu$ on $\cv$:
\[
(\grad \bu)_\cv=\frac 1 {\mcv} \int_\cv \grad \bu(x) \dx.
\]
Thanks to the regularity of $\bu$ (and the fact that $\bu=0$ on the boundary
of $\O$), there exists $\ctel{C1}$, only depending on $\bu$ (indeed, $\cter{C1}$ only depends
on the $L^\infty$-norm of the second derivatives of $\bu$),
such that, for all $\edge =K|L \in \edgesint$,
\be
\vert e_\sigma \vert  \le \cter{C1} h_\disc , \hbox{ with } e_\sigma =
(\grad \bu)_\cv \cdot \n_{\cv,\edge} - \frac {\bu(x_L)-\bu(x_K)} {d_\edge},
\label{keyrega}
\ee
and, for all $\edge \in \edgescvext$,
\be
\vert e_\sigma \vert  \le \cter{C1} h_\disc , \hbox{ with } e_\sigma =
(\grad \bu)_\cv \cdot\n_{\cv,\edge} - \frac {-\bu(x_K)} {d_{\cv,\edge}}.
\label{keyregb}
\ee
Thanks to \refe{cons}, \refe{keyrega} and \refe{keyregb}, we get, for all $\cv \in \mesh$:
\[
\dsp \mcv  (\grad_\disc P_\disc\bu)_\cv  = \sum_{\edge \in \edgescv}  \medge (x_{\edge}-x_\cv) (\grad \bu)_\cv 
\cdot\n_{\cv,\edge} + R_\cv,
\]
with $\dsp R_\cv = -\sum_{\edge \in \edgescv} e_\sigma   \medge {\rm d}(x_\edge,x_\cv)$.
Applying \refe{chxs} gives 
\be
\mcv (\grad_\disc P_\disc\bu)_\cv =   \mcv (\grad \bu)_\cv  + R_\cv.
\label{majojo}\ee
Using the inequalities \refe{keyrega} and \refe{keyregb}, we have
\be
\ba
\dsp \vert R_\cv \vert  \le  
\dsp \frac {\cter{C1}} {\theta} h_\disc  \sum_{\edge \in \edgescv}  \medge d_{\cv,\edge}=
 \frac {d\ \cter{C1}} {\theta} h_\disc  \mcv.
\ea
\label{estrk}
\ee
Then, from \refe{majojo} and \refe{estrk}, we obtain
\be
\ba
\dsp \sum_{\cv \in \mesh} \vert (\grad_\disc P_\disc\bu)_\cv - (\grad \bu)_\cv \vert ^2 \mcv \le
\\
\dsp \hfill  \sum_{\cv \in \mesh}
\left( \frac {d\ \cter{C1}} {\theta}\right)^2 h_\disc ^2 \mcv=\meas(\O) \left( \frac {d\ \cter{C1}} 
{\theta}\right)^2 h_\disc ^2.
\ea
\label{estrkb}
\ee

\medskip

In order to conclude, we remark that, thanks to the regularity of $\bu$, there exists $\ctel{C2}$,
only depending on $\bu$ (here also, $\cter{C2}$ only depends on the $L^\infty$-norm
of the second derivatives of $\bu$), such that:
\be
\sum_{\cv \in \mesh} \int_\cv \vert \grad \bu(x)-(\grad \bu)_\cv \vert^2 \dx \le \cter{C2}\ h_\disc ^2.
\label{keyregc}
\ee
Then, using \refe{estrkb} and \refe{keyregc}, we get the existence of $\cter{cstcg}$, only depending on
$\O$, $\theta$ and $\bu$, such that \refe{pdumu} holds.
\end{proof}

\begin{remark}[Choice of the points $x_K$ and $x_\edge$] \label{below} Note that in the  
proof of Lemma \ref{wkcvg}, one is free to choose any point  lying on $\edgecvcvv$
instead of  
$x_{K|L}$ in the definition of the coefficients $A_{K,L}$. 
 However, we   need this choice in the proof of 
  the strong consistency of the discrete gradient (Lemma \ref{cg}). 
Conversely, in the proof of Lemma \ref{cg}, we could take  any point of 
$K$ instead of $x_K$ in the definition of $A_{K,L}$. However, the choice of $x_K$
is crucial in the proof of Lemma \ref{wkcvg}: when comparing the terms
$\termr{A2}$ and $\termr{B}$, one needs the property of consistency of the normal flux, which 
follows from the fact that 
$\n_{K,L} = \frac {x_L - x_K} {\dkl}.$
\end{remark}


\begin{lemma}\label{convgrad}~{\bf (A sufficient condition for the strong convergence 
of the discrete gradient)}

Let $\O$ be an open bounded connected polygonal subset of $\R^d$,
$d\in\N^\star$, let $\theta >0$ and
let  $\disc$ be an admissible finite volume discretizations in the sense of  Definition
\ref{adisc}, such that $\theta_\disc \ge\theta$. Assume that there exists a function
$u_\disc\in H_\disc$ and a function $\bu\in H^1_0(\O)$ such that $u_\disc$ tends to $\bu$ in $L^2(\O)$
as $h_\disc $ tend to 0.
Assume also that there exists a function
$\alpha\in L^\infty(\O)$ and $\alpha_0 >0$ such that $\alpha(x) \ge \alpha_0$ for
a.e. $x\in\O$ and
$ [u_\disc,u_\disc]_{\disc,\alpha}$ tends to
$\int_\O  \alpha(x) \grad  \bu(x)^2 \dx$ as $h_\disc $ tends to 0. 
Then $\grad_{\disc}u_\disc$ tends to $\grad\bu$ in $L^2(\O)^d$ as $h_\disc $ tends to 0.
\end{lemma}
\begin{proof}
Let $\varphi\in {\rm C}^\infty_c(\O)$ be given (this function is devoted to approximate $\bu$ in $H^1_0(\O)$). 
Thanks to the Cauchy-Schwarz inequality, we have
\[
\int_\O (\grad_{\disc}u_\disc(x) - \grad \bu(x))^2 \dx \le 3 \
(\terml{f2}^\disc+\terml{f6}^\disc+\terml{f5})
\]
with
\[
\termr{f2}^\disc = \int_\O (\grad_{\disc}u_\disc(x) - 
\grad_{\disc}P_\disc\varphi(x))^2 \dx,
\]
\[
\termr{f6}^\disc =
\int_\O (\grad_{\disc}P_\disc\varphi(x) - \grad\varphi(x))^2 \dx,
\]
and
\[
\termr{f5} =
\int_\O (\grad\varphi(x) - \grad \bu(x))^2 \dx.
\]
We have, thanks to Lemma \ref{cg},
\be
\lim_{h_\disc \to 0} \termr{f6}^\disc = 0.
\label{vvv41}\ee
Thanks to Lemma \ref{boundnorm}, we have
\[
\int_\O (\grad_{\disc}v(x))^2 \dx \le \cter{C}^{\ 2} [v,v]_{\disc,1} \le 
\frac {\cter{C}^{\ 2}} {\alpha_0} [v,v]_{\disc,\alpha},\ \forall v\in H_{\disc}.
\]
We thus get, setting $v = u_\disc - P_\disc\varphi$ in the above inequality, that
\[
\termr{f2}^\disc \le \frac {\cter{C}^{\ 2}} {\alpha_0} ([u_\disc,u_\disc]_{\disc,\alpha}
- 2 [u_\disc,P_\disc\varphi]_{\disc,\alpha} + 
[P_\disc\varphi,P_\disc\varphi]_{\disc,\alpha}).
\]
We have, applying twice Lemma \ref{discpct}, that
\be
\lim_{h_\disc \to 0} [u_\disc,P_\disc\varphi]_{\disc,\alpha} = 
\int_\O \alpha(x) \grad \bu(x)\cdot\grad\varphi(x) \dx
\label{vvv31}\ee
and
\be
\lim_{h_\disc \to 0} [P_\disc\varphi,P_\disc\varphi]_{\disc,\alpha} = \int_\O \alpha(x) \grad\varphi(x)^2 \dx.
\label{vvv32}\ee
Under the hypotheses of the lemma, we then get that
\[
\limsup_{h_\disc \to 0}\termr{f2}^\disc \le \frac {\cter{C}^{\ 2}} {\alpha_0} 
\int_\O \alpha(x) (\grad \bu(x) - \grad\varphi(x))^2 \dx.
\]
We then get, gathering the above results, setting 
$\ctel{ggg} = \frac {\cter{C}^{\ 2}} {\alpha_0} {\rm ess}\sup_{x\in\O}\alpha(x) + 1$, that 
\[
\int_\O (\grad_{\disc}u_\disc(x) - \grad \bu(x))^2 \dx\le 
\cter{ggg} \int_\O (\grad\varphi(x) - \grad\bu(x))^2 \dx + 
\terml{f7}^\disc,
\]
with
\be
\lim_{h_\disc \to 0} \termr{f7}^\disc = 0.
\label{vvv42}\ee
Let $\varepsilon >0$. We can choose $\varphi$ such that
$\int_\O (\grad\varphi(x) - \grad\bu(x))^2 \dx \le \varepsilon$, and we can 
then
choose $h_\disc $ such that
$\termr{f7}^\disc\le \varepsilon$. This completes the proof that
\be
\lim_{h_\disc \to 0} \int_\O (\grad_{\disc}u_\disc(x) - \grad \bu(x))^2 \dx = 
0.
\label{vvv5}\ee
\end{proof}

\begin{remark}\label{rt}
Thanks to Lemma \ref{convgrad}, we get the strong convergence of the discrete
gradient in the case of the classical finite volume scheme for an isotropic
problem. Note that in the above proof, we did not use the weak convergence of
the discrete gradient, and therefore any point of $K$ can be taken instead
of $x_K$ in the definition of the coefficients $A_{K,L}$. We thus find that
the average value in $K$ of the gradient defined in \cite{convgrad} is also 
strongly convergent (the average of this gradient, defined by the generalized  
Raviart-Thomas basis 
functions, is obtained by replacing $x_K$ by the barycenter of $K$ in the definition of
$A_{K,L}$). Note that the drawback of the generalization of the Raviart-Thomas basis
was the difficulty for computing approximate values of the gradients. This drawback
no longer exists for an averaged gradient. Nevertheless, the properties
of convergence of the finite volume method  shown here for non isotropic problems 
are only
proven for the choice \refe{defaint} in the definition of  $A_{K,L}$, and  not for the 
Raviart-Thomas basis. 
\end{remark}

\section{Application to Problem \refe{ellgen}}\label{seccvglin}
 
\subsection{The finite volume scheme} \label{secfv}

Under hypotheses \refe{hypomega}-\refe{hypfg}, let $\disc$ be an admissible
discretization of $\O$ in the sense of Definition \ref{adisc}. 
The finite volume approximation to Problem \refe{ellgen} is given as the 
solution of the following equation:
\be\left\{\ba
u_\disc \in H_\disc, \\
 \int_\O  (\Lambda(x) - \alpha(x){\rm I}_d) \grad_{\disc}u_\disc(x) \cdot  \grad_{\disc}v(x) \dx 
+ [u_\disc,v]_{\disc,\alpha} = \int_\O f(x) v(x) \dx,
\ \forall v \in H_\disc,
\ea\right.\label{schvf}\ee
denoting by ${\rm I}_d$ the identity application of $\R^d$.
The existence and the uniqueness of the solution $u_\disc$ to \refe{schvf}
will be stated in Lemma \ref{exun}. 
Note that in this formulation, we use the discrete gradient on part of the 
the operator only, while on a homogeneous part, we write the usual 
cell centered scheme. This needs to be done in order to obtain the stability of 
the scheme, that is some {\it a priori} estimate on the discrete solution. 
If we take $\alpha = 0$ in \refe{schvf}, we are no longer able to prove 
the discrete $H^1$ estimate \refe{estimu} below.
Taking for $v$ the characteristic function of a control volume $K$ in \refe{schvf},
 we may 
note that Equation \refe{schvf} is equivalent to finding
the values $(u_K)_{K\in\mesh}$ (we again denote
$u_K$ instead of $(u_\disc)_K$),
solution of the following system of equations:

\be\begin{array}{ll}\dsp
\dsp
\sum_{L\in\NN_K} F_{KL} + \sum_{\edge\in\edgescvext} F_{K\edge}&=
 \dsp\int_K f(x) \dx, \ \ \dsp\forall K\in\mesh,
\end{array} 
\label{schvfdu}\ee

 where
\be\label{defflu1}
F_{KL} =  \tkl \alpha_\kl(u_K - u_L) + \left(\ba\Lambda_L
A_{LK}\cdot \grad_{\disc}u_L 
-  \Lambda_K  A_{KL}\cdot \grad_{\disc}u_K \ea\right)
 \ \forall K|L\in\edgesint,
\ee
and
\be\label{defflus1}
F_{K\edge} = \tau_{K\edge} \alpha_\edge u_K +  
\Lambda_K  A_{K\edge}\cdot \grad_{\disc}u_K
 \ \forall \edge\in\edgescvext.
\ee
In \refe{defflu1} and \refe{defflus1}, the matrices $(\Lambda_K)_{K\in\mesh}$ are defined by: 
\be
\Lambda_K  = \frac {1} {\mcv} 
\int_{K} (\Lambda(x) - \alpha(x) {\rm I}_d)  \dx.
\label{lambdaK}\ee
On can then complete the discrete expressions of $F_{KL}$ and $F_{K\edge}$ using
Definition \ref{espdis} for $A_{KL}$ $A_{K\edge}$,  and $\grad_{\disc}u_K$ for all
$K\in\mesh$, $L\in \NN_K$ and $\edge\in\edgescv$.


This is indeed a finite volume scheme, since
\[
F_{KL} = - F_{LK}, \ \forall K|L\in\edgesint.
\]

The existence of a solution to \refe{schvf} will be proven below.

\subsection{Discrete $H^1(\O)$ estimate}

We now prove the following estimate:

\begin{lemma}\label{estl1l2}
{\bf [Discrete $H^1$ estimate]}
Under hypotheses \refe{hypomega}-\refe{hypfg}, let $\disc$ be an admissible
discretization of $\O$ in the sense of Definition \ref{adisc}. 
Let $u\in H_\disc$ be a solution
to \refe{schvf}.
 Then the following inequalities hold:

\be
 \alpha_0 \Vert u\Vert_\disc \le {\rm diam}(\O) \Vert f \Vert_{(L^2(\O))^2},
\label{estimu}\ee

\end{lemma}

\begin{proof} We apply \refe{schvf} setting $v = u$ . 
We get
\[
\int_\O  (\Lambda(x) - \alpha(x) {\rm I}_d) \grad_{\disc}u(x) \cdot  \grad_{\disc}u(x) \dx 
+ [u,u]_{\disc,\alpha}= \int_\O f(x) u(x) \dx,
\]
which implies
\[
\alpha_0 [u,u]_{\disc} \le \int_\O f(x) u(x) \dx.
\]
Then the conclusion follows from the discrete Poincar\'e inequality \refe{poindis}.
\end{proof}

We can now state the existence and the uniqueness of a discrete solution
to \refe{schvf}. 
\begin{cor}\label{exun}{\bf [Existence and uniqueness of a solution to the 
finite volume scheme]}
Under hypotheses \refe{hypomega}-\refe{hypfg}, let $\disc$ be an admissible
discretization of $\O$ in the sense of Definition \ref{adisc}. Then there exists a unique
$u_\disc$ solution to \refe{schvf}.
\end{cor}
\begin{proof}
System \refe{schvf} is a linear system. Assume that $f=0$. From the discrete
Poincar\'e inequality \refe{poindis}, we get that $u = 0$. 
This proves that the linear system \refe{schvf} is invertible.
\end{proof}

\subsection{Convergence}

We have the following result, which states the convergence of the scheme
\refe{schvf}.
\begin{theo}\label{cvgce}
{\bf [Convergence of the finite volume scheme]}
Under hypotheses \refe{hypomega}-\refe{hypfg}, let $\theta>0$.  Let $\disc$
be an admissible
discretization of $\O$ in the sense of Definition \ref{adisc}, 
such that $\theta_\disc \ge \theta$.
Let $u_{\disc}\in H_{\disc}(\O)$ be the solution
to \refe{schvf}.
Then 
\begin{itemize}
\item $u_{\disc}$
converges in $L^2(\O)$ to $\bu$, weak solution of Problem \refe{ellgen}
in the sense of Definition \ref{weaksol},
\item  the discrete gradient $\grad_{\disc}u_{\disc}$ converges in
$L^2(\O)^d$ to $\grad \bu$,
\end{itemize}
 as $h_\disc $ tends to $0$.
\end{theo}

\begin{proof}
We consider a sequence of admissible discretizations
 $(\disc_n)_{n\in\N}$ such that $h_\discn $ tend to $0$ as $n\to\infty$
and  $\theta_\discn \ge \theta$
for all $n\in\N$. Thanks to Lemma \ref{estl1l2}, we can apply the compactness
result \refe{discpct}, which gives the existence of
a subsequence (again denoted $(\disc_n)_{n\in\N}$), and of 
$\bu\in H^1_0(\O)$ such that $u_{\disc_n}$ (given by \refe{schvf}
with $\disc = \disc_n$) tends to $\bu$
in $L^2(\O)$ as $n\tends\infty$. 
Let $\varphi\in {\rm C}^\infty_c(\O)$ be given,
we choose $v = P_{\disc_n}\varphi$ as test function in \refe{schvf}. We obtain
\be
 \int_\O  (\Lambda(x) - \alpha(x){\rm I}_d) \grad_{{\disc_n}}u_{\disc_n}(x) 
\cdot  \grad_{{\disc_n}}P_{\disc_n}\varphi(x) \dx 
+ [u_{\disc_n},P_{\disc_n}\varphi]_{{\disc_n},\alpha} = \int_\O f(x) P_{\disc_n}\varphi(x) \dx.
\label{schvfphi}\ee
We let $n\to\infty$ in \refe{schvfphi}. 
Thanks to Lemma \ref{wkcvg} and
Lemma \ref{cg} (which provide a weak/strong convergence result), we get that
\[
\lim_{n\to\infty}
\int_\O  (\Lambda(x) - \alpha(x) {\rm I}_d) \grad_{{\disc_n}}u_{{\disc_n}}(x) \cdot  
\grad_{{\disc_n}}P_{\disc_n}\varphi(x) \dx = 
\int_\O (\Lambda(x) - \alpha(x) {\rm I}_d) \grad  \bu(x)\cdot  \grad  \varphi(x)\dx.
\]
Using Lemma \ref{discpct}, we get that
\[
\lim_{n\to\infty} [u_{\disc_n},P_{\disc_n}\varphi]_{{\disc_n},\alpha} = 
\int_\O \alpha(x)\grad  \bu(x)\cdot  \grad  \varphi(x)\dx.
\]
Since it is easy to see that
\[
\lim_{n\to\infty} \int_\O f(x) P_{\disc_n}\varphi(x) \dx = \int_\O f(x) \varphi(x) \dx,
\]
we thus get that any limit $\bu$ of a subsequence of solutions
satisfies \refe{ellgenf} with $v = \varphi$. A classical density argument and
the uniqueness of the solution to \refe{ellgenf} permit to conclude to 
the convergence in $L^2(\O)$ of $u_{\disc}$ to $\bu$, weak solution of the problem
in the sense of Definition \ref{weaksol}, as $h_\disc $ tends  to $0$, thanks to the 
fact that  $\theta_\disc \ge \theta$.
Let us now prove the strong convergence of $\grad_{\disc}u_{\disc}$ to $\grad \bu$.
We have, using \refe{schvf} with $v = u_\disc$,
\be
\int_\O  (\Lambda(x) - \alpha(x) {\rm I}_d) \grad_{\disc}u_{\disc}(x) \cdot  
\grad_{\disc}u_{\disc}(x) \dx 
 = \int_\O f(x) u_{\disc}(x) \dx - [u_{\disc},u_{\disc}]_{\disc,\alpha}.
\label{vvv0}\ee
Thanks to Lemma \ref{discpct}, we have
\[
\int_\O \alpha(x)\grad  \bu(x)^2 \dx \le \liminf_{h_\disc \to 0} 
[u_{\disc},u_{\disc}]_{\disc,\alpha},
\]
and therefore, passing to the limit in \refe{vvv0}, we get that
\[
 \limsup_{h_\disc \to 0}\int_\O  (\Lambda(x) - \alpha(x) {\rm I}_d) 
\grad_{\disc}u_{\disc}(x) \cdot  \grad_{\disc}u_{\disc}(x) \dx 
 \le \int_\O f(x) u_{\disc}(x) \dx - \int_\O \alpha(x)\grad  \bu(x)^2 \dx.
\]
We then have, letting $v=\bu$ in \refe{ellgenf},
\be
\int_\O (\Lambda(x) - \alpha(x) {\rm I}_d) \grad  \bu(x)\cdot\grad  \bu(x) \dx = 
\int_\O f(x)  \bu(x) \dx - \int_\O \alpha(x)\grad  \bu(x)^2 \dx.
\label{vvv3}\ee
This leads to
\[
 \limsup_{h_\disc \to 0}\int_\O  (\Lambda(x) - \alpha(x) {\rm I}_d) 
\grad_{\disc}u_{\disc}(x) \cdot  \grad_{\disc}u_{\disc}(x) \dx 
 \le \int_\O (\Lambda(x) - \alpha(x) {\rm I}_d) \grad  \bu(x)\cdot\grad  \bu(x) \dx.
\]
Using Lemma \ref{wkcvg}, which states
the weak convergence of the gradient  $\grad_{\disc}u_{\disc}$ to
$\grad \bu$, we get that
\[
\int_\O (\Lambda(x) - \alpha(x) {\rm I}_d) \grad  \bu(x)\cdot\grad  \bu(x) \dx 
\le  \liminf_{h_\disc \to 0}
\int_\O  (\Lambda(x) - \alpha(x) {\rm I}_d) \grad_{\disc}u_{\disc}(x) \cdot  
\grad_{\disc}u_{\disc}(x) \dx.
\]
The above inequalities yield
\be
\lim_{h_\disc \to 0}
\int_\O  (\Lambda(x) - \alpha(x){\rm I}_d) \grad_{\disc}u_{\disc}(x) \cdot  
\grad_{\disc}u_{\disc}(x) \dx = 
\int_\O (\Lambda(x) - \alpha(x) {\rm I}_d) \grad  \bu(x)\cdot  \grad  \bu(x)\dx.
\label{vvv}\ee
{From} \refe{vvv0}, \refe{vvv3} and \refe{vvv}, we thus obtain that
\[
\lim_{h_\disc \to 0} [u_{\disc},u_{\disc}]_{\disc,\alpha} =
\int_\O  \alpha(x) \grad  \bu(x)^2 \dx,
\]
Therefore we can apply Lemma \ref{convgrad}. This completes 
the proof of the strong convergence of the 
discrete gradient.
\end{proof}

\section{Error estimate}\label{errest}

We now give an error estimate, assuming first that the solution of \refe{ellgenf}
is in ${\rm C}^2(\overline \Omega)$. In Theorem
\ref{eed}, we will consider the weaker hypothesis that the solution of \refe{ellgenf}
is only in $H^2(\Omega)$ under the assumption $d\le 3$.

\begin{theo}[${\rm C}^2$ error estimate]\label{eeu}
Assume hypotheses  \refe{hypomega}-\refe{hypfg} and that $\Lambda$ and $\alpha$ are
of class ${\rm C}^1$ on $\overline \O$.
Let  $\disc$ be an admissible finite volume discretization 
(in the sense of  Definition
\ref{adisc}).
Let $\theta\in(0,\theta_\disc ]$, where $\theta_\disc$ is defined by
\refe{reguld}.
Let $u_\disc \in H_\disc$ be the solution of \refe{schvf} and $\bu \in H^1_0(\Omega)$
be the solution of \refe{ellgenf}. We assume that $\bu \in {\rm C}^2(\overline \Omega)$.

Let us first assume that
\be 
\forall \edge \in \edgesext,
 \int_\edge \Lambda(x)
 \n_{\partial \Omega}(x) \cdot (x_\sigma - z_\sigma)  {\rm d}\gamma(x)= 0,
 \label{hypregee}\ee 
where $ \n_{\partial \Omega}(x)$ is the unit normal vector to $\partial \Omega$
at point $x$, outward to $\Omega$. 

Then, there exists $\ctel{CC2}$ only depending on $\Omega$, $\theta$,
$\alpha_0$, $\alpha$, $\beta$, $\Lambda$
and $\Vert \bu \Vert_{{\rm C}^2(\Omega)}$, such that:

\be
\Vert u_\disc - P_\disc \bu \Vert_\disc \le \cter{CC2} h_\disc ,
\label{eeuu}
\ee

\be
\Vert u_\disc - \bu \Vert_{L^2(\O)} \le \cter{CC2} h_\disc ,
\label{eeud}
\ee
and
\be
\Vert \grad_\disc u_\disc - \grad \bu \Vert_{L^2(\O)^d} \le \cter{CC2} h_\disc .
\label{eeut}
\ee

Let us then assume that  \refe{hypregee} no longer holds, then 
there exists $\ctel{cc22}$, only depending on $\Omega$, $\theta$,
$\alpha$, $\beta$, $\Lambda$
and $\Vert \bu \Vert_{H^2(\Omega)}$, such that \refe{eedu}, \refe{eedd},
\refe{eedt} hold 
with $\cter{cc22} \sqrt h_\disc$ instead of  $\cter{CC2}  h_\disc$. 

\end{theo}

\begin{remark}\label{estred}
Let us give some sufficient (and practical) 
conditions  for \refe{hypregee} to hold ~:
 
\begin{itemize}

\item If the normal vector to $\dr \O$ is an eigenvector of $\Lambda(x)$
for a.e. $x \in \dr \O$, then \refe{hypregee} holds. 
 Since this property is  always satisfied in the isotropic case, 
 the   error estimate on the gradient  \refe{eeut} holds for the 
classical cell centered scheme, for any admissible mesh.

\item  If  for all $\edge \in \edgesext$ with $\edge \in
\edgescv$, the
barycenter $x_\edge$ of $\edge$ is equal to the orthogonal projection $z_\edge$
of $x_\cv$ on $\edge$, then  \refe{hypregee} holds. 
This hypothesis is easy to ensure on 
rectangular and triangular meshes. 
 
\end{itemize}

Note also that one could   replace   \refe{hypregee}  
by $\vert z_\edge - x_\edge \vert \le \frac 1 {\theta} \diam(\cv) (h_\disc )^\half$
 for all $\edge
\in \edgesext$.


\end{remark}

\begin{proof}
In the proof, we denote by $C_i$ ($i \in \N$), various quantities only depending on
$\Omega$, $\theta$, $\alpha_0$, $\alpha$, $\beta$, $\Lambda$
and $\Vert \bu \Vert_{{\rm C}^2(\Omega)}$.

\smallskip

{\bf Step 1.} Let $v \in H_\disc$. We first perform a computation of a consistency error, namely a bound for
$\vert \terml{EC2}(v) \vert $ where $\termr{EC2}(v)$ is defined by:
\be 
 \int_\O  (\Lambda(x) - \alpha(x) {\rm I}_d) \grad_{\disc} P_\disc \bu (x) \cdot  \grad_{\disc}v(x) \dx 
+ [P_\disc \bu,v]_{\disc,\alpha} = \int_\O f(x) v(x) \dx +\termr{EC2}(v).
\label{ce}
\ee
We first consider the second term of the left hand side of \refe{ce}.
Using classical consistency error (also used in the proof
of Lemma \ref{discpct}), one has:
\be
[P_\disc \bu,v]_{\disc,\alpha} = -\int_\Omega  \div (\alpha \grad  \bu) (x) v(x) \dx + \terml{E1C2}(v),
\label{estu}
\ee
with
\[
\vert \termr{E1C2}(v) \vert \le \sum_{\edge \in \edges}  \medge \vert R_\edge \vert  
\delta_\edge v ,
\]
where $\delta_\edge v= |v_K-v_L|$ if $\edge=K|L$ is an interior edge,
$\delta_\edge v=|v_K|$ is $\edge\in\edgesext$ 
and $\vert R_\edge \vert \le \ctel{C1C2} h_\disc $. Using
the Cauchy-Schwarz inequality, this leads to:
\be
\vert \termr{E1C2}(v) \vert \le \ctel{C2C2} h_\disc  \Vert v \Vert_\disc.
\label{estuu}
\ee

\smallskip

We now consider the first term of the left hand side of \refe{ce}. We have
\be
\int_\Omega (\Lambda(x)-\alpha(x){\rm I}_d) \grad_\disc  P_\disc\bu(x) \cdot \grad_\disc v(x) \dx=
\terml{TC2}(v) + \terml{E2C2}(v),
\label{estd}
\ee
with
\[
\termr{TC2}(v) = \int_\Omega (\Lambda(x)-\alpha(x){\rm I}_d) \grad \bu(x) \cdot \grad_\disc v(x) \dx
\]
and
\[
\vert \termr{E2C2}(v) \vert  \le \ctel{C3C2} 
\Vert \grad_\disc P_\disc \bu - \grad \bu \Vert_{L^2(\O)^d} \Vert \grad_\disc v \Vert_{L^2(\O)^d}.
\]
Using Lemma \ref{cg} and  Lemma \ref{boundnorm}, we obtain
\be
\vert \termr{E2C2}(v) \vert  \le \ctel{C4C2} h_\disc  \Vert v \Vert_\disc.
\label{estdd}
\ee

We now compute $\termr{TC2}(v)$.
For $\cv \in \mesh$ and $\edge \in \edges$, let $\mu_\cv$ and $\mu_\edge$ respectively be
the mean values of $(\Lambda(x)-\alpha(x){\rm I}_d)
\grad \bu$ on $\cv$ and $\edge$:
\[
\mu_\cv= \frac 1 {\mcv} \int_\cv (\Lambda(x)-\alpha(x){\rm I}_d)
\grad \bu (x) \dx, \; \;
\mu_\edge= \frac 1 {\medge} \int_\edge (\Lambda(x)-\alpha(x){\rm I}_d
\grad \bu (x) {\rm d}\gamma(x).
\]
The regularity of $\bu$, $\Lambda$ and $\alpha$ gives, for all $\cv \in \mesh$
and all $\edge \in \edgescv$ (recall that $\vert \cdot \vert$ denotes the Euclidean norm in $\R^d$):
\be
\vert \mu_\cv-\mu_\edge \vert \le \ctel{C5C2} h_\disc .
\label{mumu}
\ee
Indeed, $\cter{C5C2}$ only depends on the $L^\infty$-norms of $\Lambda$, $\alpha$
and $\grad \bu$ and on the $L^\infty$-norms of the derivatives of $\Lambda$, $\alpha$
and $\grad \bu$.

\smallskip

We now use \refe{mumu} in order to give a bound  of $\termr{TC2}(v)$
as a function of $h_\disc $. 
Indeed, the definition of $\grad_\disc v$ leads to:
\[
\ba
\dsp \termr{TC2}(v)= \sum_{K \in \mesh} \mu_\cv \cdot \mcv (\grad_\disc v)_K=\\ \dsp
\sum_{K \in \mesh} 
\left(\sum_{L\in\NN_K}  \mu_\cv \cdot A_{K,L} \ (v_L-v_K)-
\sum_{\edge\in \edgescvext}  \mu_\cv \cdot A_{K,\edge} v_K\right)=
\\
\dsp 
\sum_{K \in \mesh} \left(\sum_{L\in\NN_K}  \mu_{K \vert L} \cdot A_{K,L} \ (v_L-v_K)-
\sum_{\edge\in \edgescvext} \mu_\edge \cdot A_{K,\edge} v_K\right) +\terml{E3C2}(v),
\ea
\]
with 
\[
\ba
\dsp \vert \termr{E3C2}(v) \vert \le \cter{C5C2} h_\disc  \sum_{K \in \mesh} \left(\sum_{L\in\NN_K}  \vert A_{K,L} \vert
\vert v_L-v_K\vert+
\sum_{\edge\in \edgescvext} \vert A_{K,\edge} \vert \vert  v_K \vert \right) \le
\\
\dsp \cter{C5C2} h_\disc  \left( \sum_{\edge=K|L \in \edgesint} 
(\vert A_{K,L} \vert + \vert A_{L,K} \vert)
\vert v_L-v_K\vert +
\sum_{K \in \mesh} \sum_{\edge\in \edgescvext} \vert A_{K,\edge} \vert \vert  v_K \vert \right).
\ea
\]

\smallskip

Since $A_{K,L} = \tkl (x_{K \vert L} - x_\cv)$ and $A_{\cv,\edge} = \tau_\edge(x_\edge-x_\cv)$, one deduces from
the preceding inequality, thanks to the definition of $\theta_\disc $ (which gives
${\rm d}(x_\edge,x_\cv) \le (d_{\cv,\edge}/\theta)$ if $\edge \in \edgescv$) and using
Cauchy-Schwarz Inequality:
\be
\vert \termr{E3C2}(v) \vert \le \ctel{C6C2} h_\disc  \Vert v \Vert_\disc.
\label{estdt}
\ee
We now remark that:
\be
\ba
\dsp \termr{TC2}(v)-\termr{E3C2}(v)=\sum_{K \in \mesh} \left(\sum_{L\in\NN_K}  \mu_{K \vert L} \cdot A_{K,L} \ (v_L-v_K)-
\sum_{\edge\in \edgescvext} \mu_\edge \cdot A_{K,\edge} v_K\right)= 
\\
\hfill
\dsp \sum_{\edge=K|L \in \edgesint} \mu_\edge \cdot (x_L-x_K) \transedge (v_L-v_K)
-\sum_{K \in \mesh} \sum_{\edge\in \edgescvext} \mu_\edge \cdot (x_\edge-x_\cv) \transedge v_\cv.
\ea
\label{ttt}
\ee

For $\edge \in \edgesint$, one has $\edge=K \vert L$ and $(x_L-x_K)=d_{\edge}\n_{K,\edge}$ where $\n_{K,\edge}$ is
the normal vector to $\edge$ exterior to $\cv$.

\smallskip

For $\edge \in \edgesext$, one has $\edge \in \edgescv$.
Thanks to the fact that under homogeneous Dirichlet boundary conditions, 
the gradient of $\bar u$ is normal to the boundary,    using Assumption 
\refe{hypregee},  we get  that
 
 $$\mu_\edge \cdot (x_\edge-x_\cv)  \transedge = 
\int_\edge (\Lambda(x)-\alpha(x){\rm I}_d 
\grad \bu (x) \cdot \n_{\partial \Omega}(x) {\rm d}\gamma(x).
$$

Then, 
one deduces from \refe{ttt}:

\begin{equation}
\termr{TC2}(v)-\termr{E3C2}(v)=-
\int_\O \div ((\Lambda-
\alpha {\rm I}_d) \grad \bu)(x) v(x) \dx.
\label{tttandt}
\end{equation}

Therefore, since $-\div( \Lambda \grad \bu)=f$,
one has \refe{ce} with $\termr{EC2}(v)=\termr{E1C2}(v)+\termr{E2C2}(v)+\termr{E3C2}(v)$. This gives,
with \refe{estuu}, \refe{estdd}, \refe{estdt}:
\be
\vert \termr{EC2}(v) \vert \le \ctel{C7C2} h_\disc  \Vert v \Vert_\disc.
\label{estf}
\ee
This concludes Step 1.

{\bf Step 2.}

Let $e_\disc=P_\disc \bu - u_\disc$ be the discrete discretization error. Using \refe{ce} and \refe{schvf} give,
for all $v \in H_\disc$:
\[
 \int_\O  (\Lambda(x) - \alpha(x) {\rm I}_d) \grad_{\disc} e_\disc(x) \cdot  \grad_{\disc}v(x) \dx 
+ [e_\disc,v]_{\disc,\alpha} = \termr{EC2}(v).
\] 
Taking $v=e_\disc$ in this formula gives, with \refe{estf}, $[e_\disc,e_\disc]_{\disc,\alpha} \le 
\cter{C7C2} h_\disc 
\Vert e_\disc \Vert_\disc$ and
then, with $\ctel{C8C2}=\cter{C7C2}/ \alpha_0$ (since $\alpha_0 \Vert e_\disc \Vert_\disc^2 \le
[e_\disc,e_\disc]_{\disc,\alpha}$):
\be
\Vert e_\disc \Vert_\disc \le \cter{C8C2} h_\disc ,
\label{eest}
\ee
which is exactly  \refe{eeuu}. 

\smallskip

Using the Discrete Poincar\'e
Estimate \refe{poindis} and the fact that  $\bu \in C(\overline \O)$, one deduces \refe{eeud} from \refe{eeuu}.

\smallskip

The last estimate, Estimate \refe{eeut}, is a direct consequence of
\refe{eest}, \refe{pdumu} and  \refe{mn}. 
This concludes the first part of the theorem, {\it i.e.} assuming 
\refe{hypregee}.

\bigskip

If $\disc$  no longer satisfies the hypothesis \refe{hypregee}, 
one has to replace \refe{tttandt} by:
\[
\termr{TC2}(v)-\termr{E3C2}(v)=-
\int_\O \div ((\Lambda-
\alpha {\rm I}_d) \grad \bu)(x) v(x) \dx + \terml{E4C2}(v),
\]
where, recalling that by $z_{\edge}$ the orthogonal projection of $x_\cv$ on
$\edge$ (see Definition~\ref{adisc}):
\[
\termr{E4C2}(v) = \sum_{K \in \mesh} \sum_{\edge\in \edgescvext} \mu_\edge \cdot (z_{\edge}-x_\edge) 
\transedge v_\cv.
\]
Thanks to the Cauchy-Schwarz inequality, we get
\[
\termr{E4C2}(v)^2 \le \sum_{K \in \mesh} \sum_{\edge\in \edgescvext} \transedge
\mu_\edge^2  (\diam(\cv))^2
 \sum_{K \in \mesh} \sum_{\edge\in \edgescvext} \transedge v_\cv^2,
\]
which leads to
\[
\termr{E4C2}(v)^2 \le \frac {h_\disc } {\theta}  \meas(\partial \Omega) \Vert \grad \bar u \Vert_\infty^2
\Vert v \Vert_\disc^2,
\]
 where  $\meas(\partial \O)$ is the $d-1$-dimensional Lebesgue measure of $\partial \O$.
This gives  \refe{estf} with $h_\disc ^\half$ instead of $h_\disc $. Following
Step 2, this  allows to conclude  the proof. 
\end{proof}

\bigskip

We now want an error estimate when the solution of \refe{ellgenf} is in $H^2(\O)$ instead of
${\rm C}^2(\overline \O)$, in the case where the
space dimension is lower or equal to 3. 
Indeed, the ${\rm C}^2$-regularity of the solution of \refe{ellgenf}
was used, in the preceding proofs,
only four times, namely to prove \refe{keyrega}, \refe{keyregb} and \refe{keyregc} in Lemma
\ref{cg} and to prove \refe{mumu} in Theorem \ref{eeu} (in fact, it is also used for the
classical consistency error \refe{estu}, but, for this term, the generalization to the case where the 
solution of \refe{ellgenf} is in $H^2(\O)$ instead of
${\rm C}^2(\overline \O)$, in the case $d\le 3$,  is already done in \cite{book}).
We will now prove similar inequalities for $\bu \in H^2(\O) \cap H^1_0(\O)$ (instead
of $\bu \in {\rm C}^2(\O)$ with $\bu=0$ on the boundary of $\O$) which will allow us
to obtain the desired error estimate.

\begin{lemma}[Consistency of the gradient, $\bu \in H^2(\O)$]\label{cgb}
Under hypothesis  \refe{hypomega}, with $d\le 3$,
let  $\disc$ be an admissible finite volume discretization in the sense of  Definition
\ref{adisc}, and let $\theta\in(0,\theta_\disc ]$.
Let $\bu \in H^2(\Omega)\cap H^1_0(\O)$.
Then, there exists $\ctel{CC3}$, only depending on $\Omega$,
$\theta$ and $\bu$, such that:
\be
\Vert \grad_\disc (P_\disc \bu) - \grad \bu \Vert_{L^2(\O)^d} \le \cter{CC3} h_\disc 
\Vert \bu \Vert_{H^2(\Omega)}.
\label{pdumubis}
\ee
(Recall that $P_\disc$ is defined in \refe{Pdisc} and $\grad_\disc$ in Definition \ref{espdis}.)
\end{lemma}

\begin{proof}

The proof follows the proof of Lemma \ref{cg} (in particular, recall that
$H^2(\O) \subset C(\overline \O)$ since $d \le 3$).
The ${\rm C}^2$-regularity was only used 
to prove  \refe{keyrega},  \refe{keyregb}, \refe{keyregc}.
We now prove similar inequalities in the case $\bu \in H^2(\O)$.

\smallskip

We begin with providing inequalities similar to  \refe{keyrega},  \refe{keyregb}.
We denote by $(\grad \bu)_\edge$ the mean value of $\grad \bu$ on $\edge$ (recall
that $(\grad \bu)_\cv$ is the mean value of $\grad \bu$ on $\cv$).
We use Inequality (9.63) of \cite{book} (in the
proof of Theorem 9.4, using the $H^2$-regularity). This inequality states the existence
of $\ctel{c0}$, only depending on $d$ and $\theta$, such that, for all $\edge =K\vert L \in \edgesint$:

\be
\vert E_\sigma \vert^2   \le \cter{c0} \frac {h_\disc ^2} {\medge d_\edge} \int_{D_\edge} \vert H(\bu)(z) \vert^2 \drm z,
\hbox{ with } E_\sigma= (\grad \bu)_\edge \cdot \n_{\cv,\edge} - \frac {\bu(x_L)-\bu(x_K)} {d_\edge},
\label{keyregabisn}
\ee
and, for all $\edge \in \edgesext$, if $\edge \in \edgescv$:
\be
\vert E_\edge \vert^2  \le
\cter{c0} \frac {h_\disc ^2} {\medge d_\edge} \int_{D_\edge} \vert H(\bu)(z) \vert^2 \drm z,
\hbox{ with } E_\edge=
(\grad \bu)_\edge \cdot \n_{\cv,\edge} - \frac {-\bu(x_K)} {d_{\cv,\edge}},
\label{keyregbbisn}
\ee
where:
\[
\vert H(\bu)(z) \vert^2 = \sum_{i,j=1}^d \vert D_iD_j\bu(z) \vert^2.
\]
We have now to compare $(\grad \bu)_\edge$ and $(\grad \bu)_\cv$. This is possible
thanks to Inequality (9.38) in Lemma 9.4 of \cite{book}. Following this result, there exists $\ctel{D}$, only
depending on $d$ and $\theta$,
such that, for all $\cv \in \mesh$, all $\edge \in \edgescv$ and all $v \in H^1(\cv)$:
\be
\ba
\dsp \left\vert \frac 1 {\mcv} \int_\cv v(x)\dx - \frac 1 {\medge} \int_\edge v(x){\rm d}\gamma(x) \right\vert^2
\le \cter{D} \frac {\diam(\cv)} {\medge} \int_\cv \vert \grad v(x) \vert^2 \dx \le
\\
\dsp \hfill 
2 \cter{D} \frac { h_\disc ^2} {\medge d_\edge} \int_\cv \vert \grad v(x) \vert^2 \dx.
\ea
\label{smumubisn}
\ee
Using \refe{smumubisn} with the derivatives of $u$,  one deduces from \refe{keyregabisn} and
\refe{keyregbbisn}, that there exists some real value $\ctel{C0}$ only depending on $d$ and $\theta$ such that

\be
\vert e_\edge \vert^2   \le \cter{C0} \frac {h_\disc ^2} {\medge d_\edge} \int_{D_\edge} \vert H(\bu)(z) \vert^2 \drm z,
\hbox{ with } e_\edge= (\grad \bu)_\cv \cdot \n_{\cv,\edge} - \frac {\bu(x_L)-\bu(x_K)} {d_\edge},
\label{keyregabis}
\ee
and, for all $\edge \in \edgesext$, if $\edge \in \edgescv$:
\be
\vert e_\edge \vert^2  \le
\cter{C0} \frac {h_\disc ^2} {\medge d_\edge} \int_{D_\edge} \vert H(\bu)(z) \vert^2 \drm z,
\hbox{ with } e_\edge=
(\grad \bu)_\cv \cdot \n_{\cv,\edge} - \frac {-\bu(x_K)} {d_{\cv,\edge}},
\label{keyregbbis}
\ee

Since $\dsp \vert R_\cv \vert \le \sum_{\sigma \in \edgescv} \frac {\medge d_{\cv,\edge}} {\theta}
\vert e_\edge \vert$ (where $R_K$ is defined in \refe{majojo}), using 
the Cauchy-Schwarz Inequality, \refe{keyregabis} and \refe{keyregbbis}
lead to the following
bound:
\[
\ba
\dsp R_\cv^2 \le \frac 1 {\theta^2} \sum_{\sigma \in \edgescv} \medge d_{\cv,\edge}
\sum_{\sigma \in \edgescv} \medge d_{\cv,\edge} e_\sigma^2 \le
\\
\dsp \hfill \frac {d \mcv} {\theta^2} \sum_{\sigma \in \edgescv} \medge d_{\cv,\edge}
\cter{C0} \frac {h_\disc ^2} {\medge d_\edge} \int_{D_\edge} \vert H(\bu)(z) \vert^2 \drm z
\ea
\]
and, since $d_{\cv,\edge} \le d_\sigma$ and $\theta_\disc  \ge \theta$:

\[
(\frac {R_K} {\mcv})^2 \mcv \le 
\frac {d\ \cter{C0}} {\theta^2} h_\disc ^2  \sum_{\sigma \in \edgescv}
\int_{D_\edge} \vert H(\bu)(z) \vert^2 \drm z.
\]

Then, \refe{estrkb} becomes:

\[
\ba
\dsp \sum_{\cv \in \mesh} \vert (\grad_\disc P_\disc \bu)_\cv - (\grad \bu)_\cv \vert^2 \mcv \le
\\
\dsp \hfill  \sum_\cv
\frac {d\ \cter{C0}} {\theta^2} h_\disc ^2
\sum_{\sigma \in \edgescv} \int_{D_\edge} \vert H(\bu)(z) \vert^2 \drm z,
\ea
\]

which gives the existence of $\ctel{C1H2}$, only depending on $d$ and $\theta$  such that:
\be
\dsp \sum_{\cv \in \mesh} \vert (\grad_\disc P_\disc \bu)_\cv - (\grad \bu)_\cv \vert^2 \mcv \le
\cter{C1H2} h_\disc ^2 \Vert \bu \Vert_{H^2(\Omega)}^2.
\label{estrkbbis}
\ee

We have now to obtain an inequality similar to \refe{keyregc} (but without using
$\bu \in {\rm C}^2(\overline \O)$). We will use here the fact that $d_{\cv,\edge} \ge \theta \diam(\cv)$
if $\edge \in \edgescv$.

\smallskip

If $\omega$ is a convex, bounded, open subset of $\R^d$, the well-known
``Mean Poincar\'e Inequality" gives, for all $v \in H^1(\omega)$:
\be
\int_\omega \vert  v(x) - m_\omega v \vert^2 \dx \le \frac 1 {m(\omega)}
d_\omega^2 m(B(0,d_\omega)) \int_\omega \vert \grad v (x) \vert^2 \dx,
\label{mpi}
\ee
where $m_\omega(v)$ is the mean value of $v$ on $\omega$, $d_\omega$ is the
diameter of $\omega$, $B(a,\delta)$
is the ball in $\R^d$ of center $a$ and radius $\delta$ and
$m(\omega)$ (resp. $m(B(a,\delta))$ is the $d$-dimensional Lebesgue
measure of $\omega$ (resp. $B(a,\delta)$).
(A discrete counterpart of \refe{mpi} is given, for instance, in \cite{book}, Lemma 10.2.)

\smallskip

Let $\cv \in \mesh$. We will  use \refe{mpi} for $\omega=\cv$. Since $d_{\cv,\edge}$ is
the distance between $x_\cv$ to $\edge$ (for $\edge \in \edgescv)$, there exists
$\edge \in \edgescv$ such that $B(x_\cv,d_{\cv,\edge}) \subset \cv$. Then, one
has $m(B(0,1))d_{\cv,\edge}^d=m(B(x_\cv,d_{\cv,\edge})) \le \mcv$ and,
using $d_{\cv,\edge} \ge \theta \diam(\cv)$, one obtains:
\be
\mcv \ge m(B(0,1)) (\theta)^d (\diam(\cv))^d.
\label{reguldd}
\ee

Taking $\omega=\cv$ in \refe{mpi},  gives, for all $\cv \in \mesh$ and all
$v \in H^1(\cv)$:
\be
\int_\cv \vert  v(x) - m_\omega v \vert^2 \dx \le \frac 1 {\theta^d}
\diam(\cv)^2  \int_\cv \vert \grad v (x) \vert^2 \dx,
\label{mpib}
\ee

Taking $v$ equal to the derivatives
of $\bu$ (which are in $H^1(\cv)$ for all $\cv \in \mesh$) in \refe{mpib}
gives the existence
of $\ctel{C2H2}$, only depending on $d$ and $\theta$, such that:

\be
\sum_{\cv \in \mesh} \int_\cv \vert \grad \bu(x)-(\grad \bu)_\cv \vert^2 \dx \le 
\cter{C2H2} h_\disc ^2
\Vert \bu \Vert_{H^2(\Omega)}^2.
\label{keyregcbis}
\ee
Then, we conclude as in Lemma \ref{cg}, using \refe{estrkbbis} and \refe{keyregcbis}, that
there exists $\cter{CC3}$ only depending on
$\O$, $\theta$ and $\bu$ such that \refe{pdumubis} holds.
\end{proof}

\begin{theo}[$H^2$ error estimate]\label{eed}
Assume hypotheses  \refe{hypomega}-\refe{hypfg} with $d\le 3$, 
and that $\Lambda$ and $\alpha$ are
of class ${\rm C}^1$ on $\overline \O$.
Let  $\disc$ be an admissible finite volume discretization in the sense of  Definition
\ref{adisc}, and let $\theta\in (0,\theta_\disc ]$. 
We assume that
that $\card(\edgescv) \le \frac 1 {\theta}$ for all $\cv \in \mesh$.
Let $u_\disc \in H_\disc$ be the solution of \refe{schvf} and $\bu \in H^1_0(\Omega)$
be the solution of \refe{ellgenf}. We assume that $\bu \in H^2(\Omega)$ (which
is necessarily true if $\O$ is convex).

Let us first assume that
Hypothesis \refe{hypregee} holds. Then, there exists 
$\ctel{CH2}$, only depending on $\Omega$, $\theta$,
$\alpha$, $\beta$, $\Lambda$
and $\Vert \bu \Vert_{H^2(\Omega)}$, such that:

\be
\Vert u_\disc - P_\disc \bu \Vert_\disc \le \cter{CH2} h_\disc ,
\label{eedu}
\ee

\be
\Vert u_\disc - \bu \Vert_{L^2(\O)} \le \cter{CH2} h_\disc ,
\label{eedd}
\ee
and
\be
\Vert \grad_\disc u_\disc - \grad \bu \Vert_{L^2(\O)^d} \le \cter{CH2} h_\disc .
\label{eedt}
\ee
(Recall that $H_\disc$, $\grad_\disc$ and $\Vert \cdot \Vert_\disc$ are defined in Definition \ref{espdis},
$P_\disc$ is defined in \refe{Pdisc}.)

Let us then assume that  \refe{hypregee} no longer holds, then 
there exists $\ctel{CH22}$, only depending on $\Omega$, $\theta$,
$\alpha$, $\beta$, $\Lambda$
and $\Vert \bu \Vert_{H^2(\Omega)}$, such that \refe{eedu}, \refe{eedd},
\refe{eedt} hold 
with $\cter{CH22} \sqrt h_\disc$ instead of  $\cter{CH2}  h_\disc$. 
\end{theo}

\begin{proof}

The proof of Theorem \ref{eed} follows the proof of Theorem \ref{eeu}. The quantities
$\cter{C1H2}$ and $\cter{C2H2}$, depending on $\theta$, are now
used to get a bound for $\termr{E1C2}(v)$  (as in [1]), and
the quantity $\cter{C4C2}$, also
depending on $\theta$ since it is obtained with \refe{pdumubis} (Lemma \ref{cgb}) instead of
\refe{pdumu} (Lemma \ref{cgb}), is used to obtain a bound for $\termr{E2C2}(v)$.

\smallskip

In order to obtain a bound for $\termr{E3C2}(v)$ (and then to conclude the proof of Theorem \ref{eed}),
we need to obtain an inequality similar to \refe{mumu} (where the ${\rm C}^2$-regularity
of $\bu$ was used), which gives a bound for the difference between
the mean values of $(\Lambda(x)-\alpha(x){\rm I}_d)
\grad \bu$ on $\cv$ and on $\edge$ if $\edge \in \edgescv$. 
Here, we will obtain a bound  for the difference between
these mean values using once again the consequence \refe{smumubisn} of
Inequality (9.38) in Lemma 9.4 of \cite{book}.
Applying \refe{smumubisn} to the derivatives of $(\Lambda-\alpha {\rm I}_d) \grad \bu$, there exists
$\ctel{C5H2}$ only depending on $\O$, $\theta$, $\Lambda$ and $\alpha$ (indeed, the
${\rm C}^1$-norms of $\Lambda$ and $\alpha$), such that, for all
$\cv \in \mesh$, all $\edge \in \edgescv$ and all $v \in H^1(\cv)$:
\be
\vert \mu_\cv-\mu_\edge \vert^2 \le \cter{C5H2} \frac {\diam(\cv)} {\medge} \Vert \bu \Vert_{H^2(\cv)}^2.
\label{mumubis}
\ee

Following the proof of Theorem \ref{eeu}, \refe{mumubis} is used to obtain
a bound for $\termr{E3C2}(v)$:

\[
\ba
\dsp  \vert \termr{E3C2}(v) \vert  \le \sum_{K \in \mesh} \left(\sum_{L\in\NN_K} 
 \vert \mu_{K\vert L} - \mu_\cv \vert \vert A_{K,L} \ (v_L-v_K) \vert  + 
\sum_{\edge\in \edgescvext} \vert \mu_\edge - \mu_K \vert  \vert A_{K,\edge} v_K \vert \right) \le
\\
\dsp \hfill \sum_{\edge=\cv \vert \cvv \in \edgesint}
\frac {\vert \mu_\edge -\mu_\cv \vert + \vert \mu_\edge -\mu_\cvv \vert} {\theta} \medge d_\edge
\frac {\delta_\edge v } {d_\edge}
+ \sum_{\edge \in \edgesext}
\frac {\vert \mu_\edge -\mu_\cv \vert} {\theta} \medge d_\edge
\frac {\delta_\edge v} {d_\edge},
\ea
\]
where, in the last term, $\cv$ is such that $\edge \in \edgescv$ and where
$\delta_\edge v =\vert v_\cv -v_\cvv \vert $ if $\edge=\cv \vert \cvv \in \edgesint$ and
$\delta_\edge v =\vert v_\cv \vert $ if $\edge= \in \edgesext \cap \edgescv$. (We
also used the fact  that $\vert A_{K,L} \vert \le \frac {\medge}
{\theta}$ and $\vert A_{\cv,\edge} \vert \le \frac {\medge} {\theta}$, thanks to $\theta_\disc  \ge
\theta$.)

\smallskip

Then, using Cauchy-Schwarz Inequality and \refe{mumubis}, one obtains:
\[
\ba
\dsp  \vert \termr{E3C2}(v) \vert  \le \Vert v \Vert_\disc \frac {\sqrt{2C_5}} {\theta} 
( \sum_{\edge=\cv \vert \cvv \in \edgesint}
d_\edge (\diam(\cv) \Vert \bu \Vert_{H^2(\cv)}^2+ \diam(\cvv) \Vert \bu \Vert_{H^2(\cvv)}^2)
\\
\hfill \dsp +  \sum_{\edge \in \edgesext} d_\edge  \diam(\cv) \Vert \bu \Vert_{H^2(\cv)}^2)^\half.
\ea
\]

Using $d_\edge \le 2 h_\disc $, $\diam(\cv) \le h_\disc $ and
the fact that $\card(\edgescv) \le \frac 1 {\theta}$ for all $\cv \in \mesh$, one deduces the
existence of $C_6$, only depending on $\O$, $\theta$, $\Lambda$ and $\alpha$, such that:

\be
\vert \termr{E3C2}(v) \vert \le C_6 h_\disc  \Vert \bu \Vert_{H^2(\O)} \Vert v \Vert_\disc.
\label{estlast}
\ee

Then, we conclude the proof of Theorem \ref{eed} exactly as in the proof
of Theorem \ref{eeu} (\refe{estlast} replaces \refe{estdt}). \end{proof}


\section{Numerical results}\label{numex}
The scheme was tried for various academic problems, for which the
analytical solution is known. 
For the Laplace equation, we compared the classical cell centered scheme to 
the new scheme, which we shall call the gradient scheme in the 
sequel. First note that in the classical 
cell centered scheme, the equation relative to a given
cell involves the neighbors of this cell, while in the gradient scheme, it
 involves the neighbors of this cell and the neighbors of the neighbors.
Hence in the case of 
a rectangular (resp. parallelipedic) mesh, the classical cell centered scheme 
is a 5 points (resp. 7 points) scheme,
while the gradient scheme is a 13 points (resp. 24 points) scheme, 
 Similarly, if one uses a triangular  (resp. tetrahedral) 
 mesh  the classical scheme is a   4 points (resp. 7 points) scheme, 
while the gradient scheme  is a   10 points (resp. at most 17 points) scheme. 
Hence the gradient scheme is more expensive in terms of time and memory,
although this is not so much, for example compared to the use
of a $Q^1$ finite element in the case of a parallelipedic mesh, 
which leads to a 27 points scheme. 

%


We   tested the gradient scheme for some real anisotropic problems, the number of cells 
varying from 100 to 6400 in the rectangular meshes case (in fact, rectangles are
squares), and from 700 to 17500 in the 
 triangular meshes case. The convergence rates have been computed by fitting a less-square regression
on the logarithmic values of the errors and of the characteristic size of the mesh. 

\medskip

The first case is an anisotropic homogenous problem with diffusion matrix 
$$\Lambda =  \left(\begin{array}{cc}
1.5 & 0.5 \\ 0.5 & 1.5 \end{array} \right).$$
The second case is a rotating permeability field, that is, the diffusion matrix
is constant in the $(r,\theta)$ coordinates and equal to 
$\Lambda_{r,\theta}
 = \dsp \left(\begin{array}{cc} 10 & .2 \\ .2 & 10 \end{array}\right).$
The exact solution is taken to be
 $u(x_1,x_2) = \frac 1 2 \ln((x_1-.5)^2 + (x_2-1.1)^2),$
 on the domain $\Omega=]0,1[\times]0,1[$. 
 The orders of convergence which were found are given Table \ref{txconvanis}.

\begin{table} 

  \begin{center}

\begin{tabular}{| l | c | c | c | c |}
\cline{2-5}
\multicolumn{1}{ c |}{} & \multicolumn{2}{| c |}{ \begin{tabular}{c} Case 1\\ 
homogeneous anisotropic\end{tabular}} 
& \multicolumn{2}{|
c |}{\begin{tabular}{c} Case 2\\heterogeneous anisotropic\end{tabular}}  \\
\cline{2-5}
\multicolumn{1}{ c |}{} & Rectangles FV 13 & {Triangles VF10} 
 & Rectangles FV 13  &  {Triangles VF10}  \\
\hline
$u$       &  2.00   &  2.0  &  2.2  &  2.0     \\
\hline
$\grad u$ &  1.00   &  1.0  &  1.4 &  1.3  \\
\hline
\end{tabular}
\caption{Rates of convergence  of FV13 and FV10 
in a homogeneous anisotropic case
and in a heterogeneous anisotropic case}\label{txconvanis}
\end{center}
\end{table} 

\begin{figure}{htb}
\begin{center}\begin{tabular}{cc}
\includegraphics[width=0.45 \linewidth]{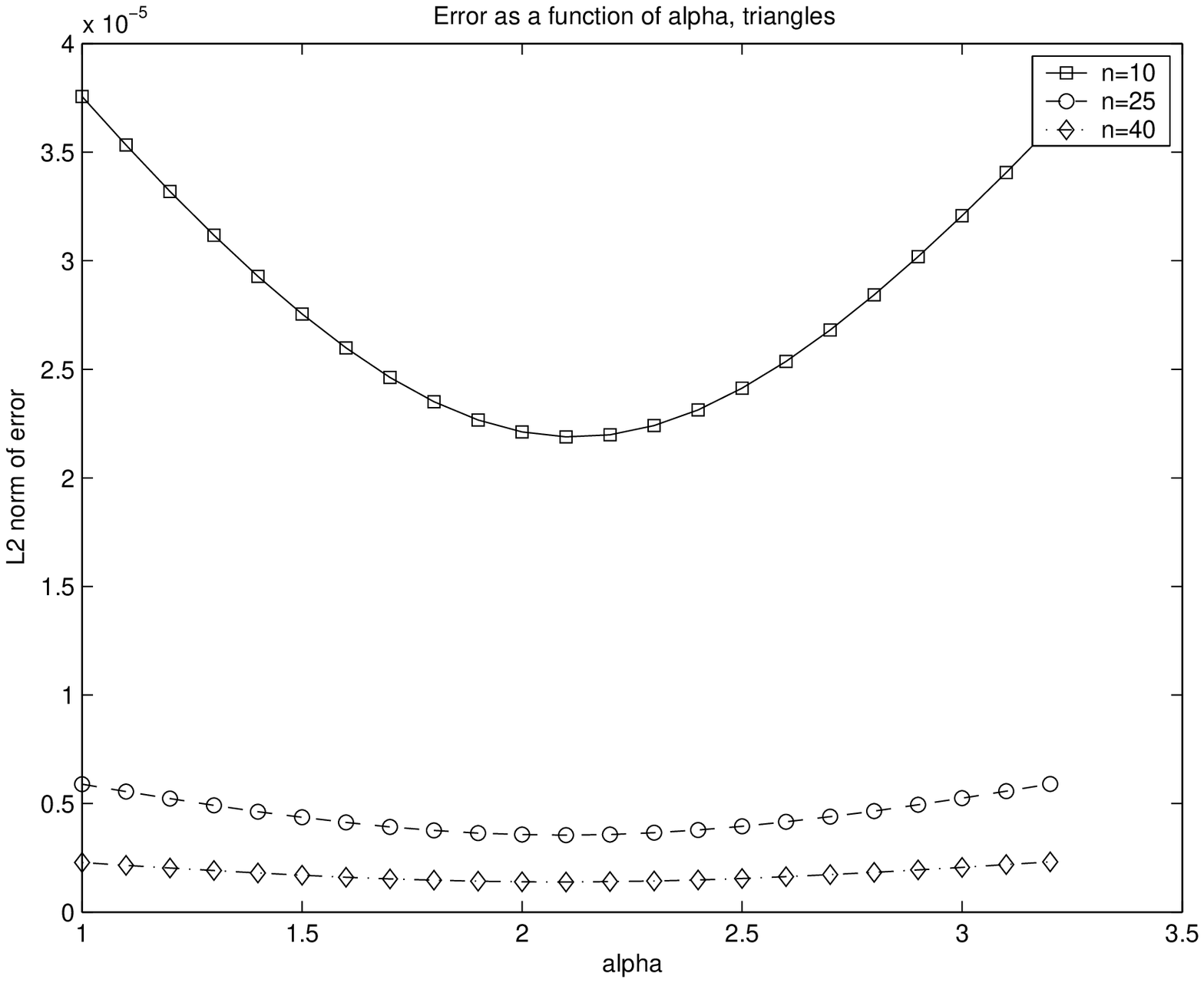}&
\includegraphics[width=0.45 \linewidth]{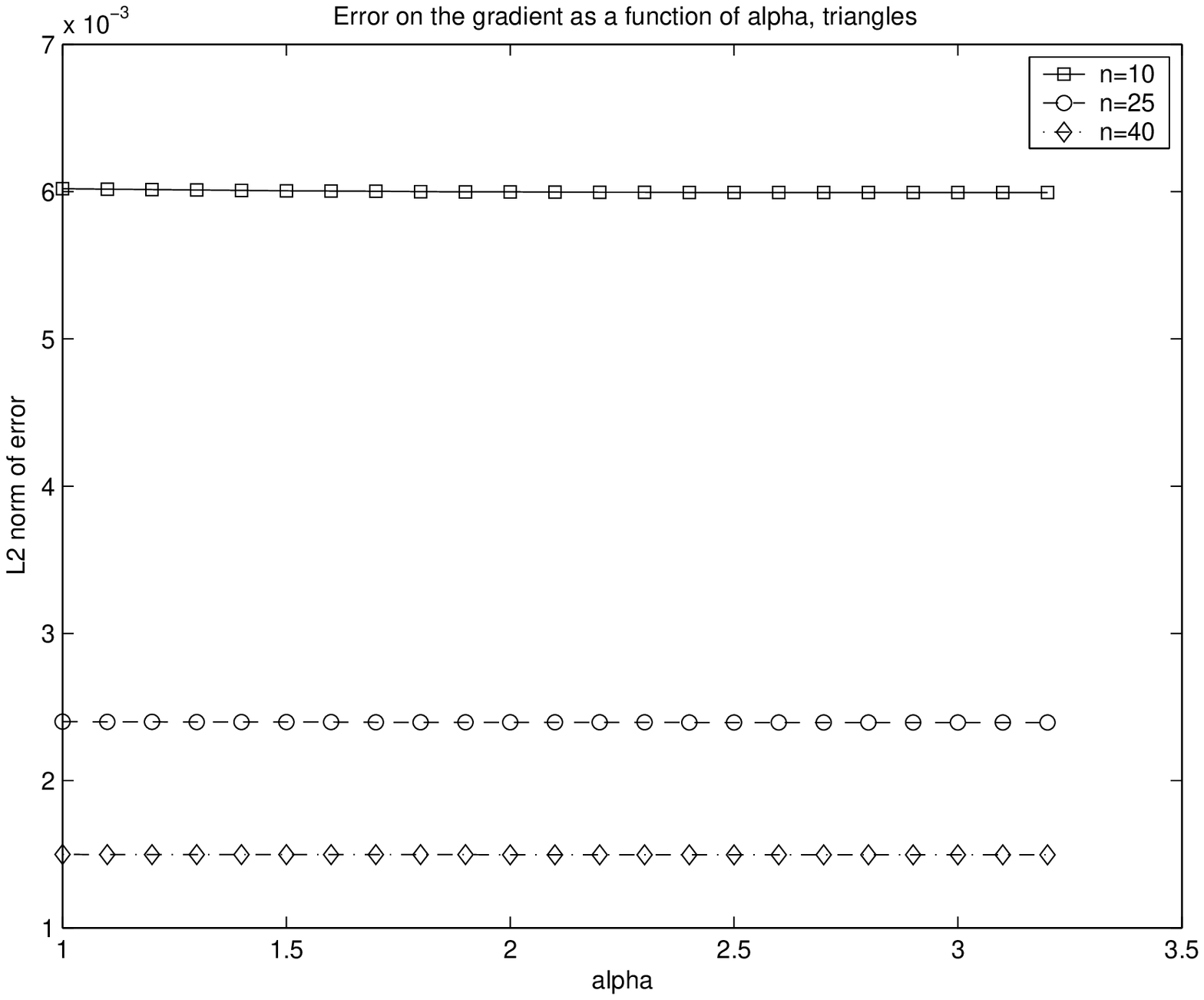}\\
\includegraphics[width=0.45 \linewidth]{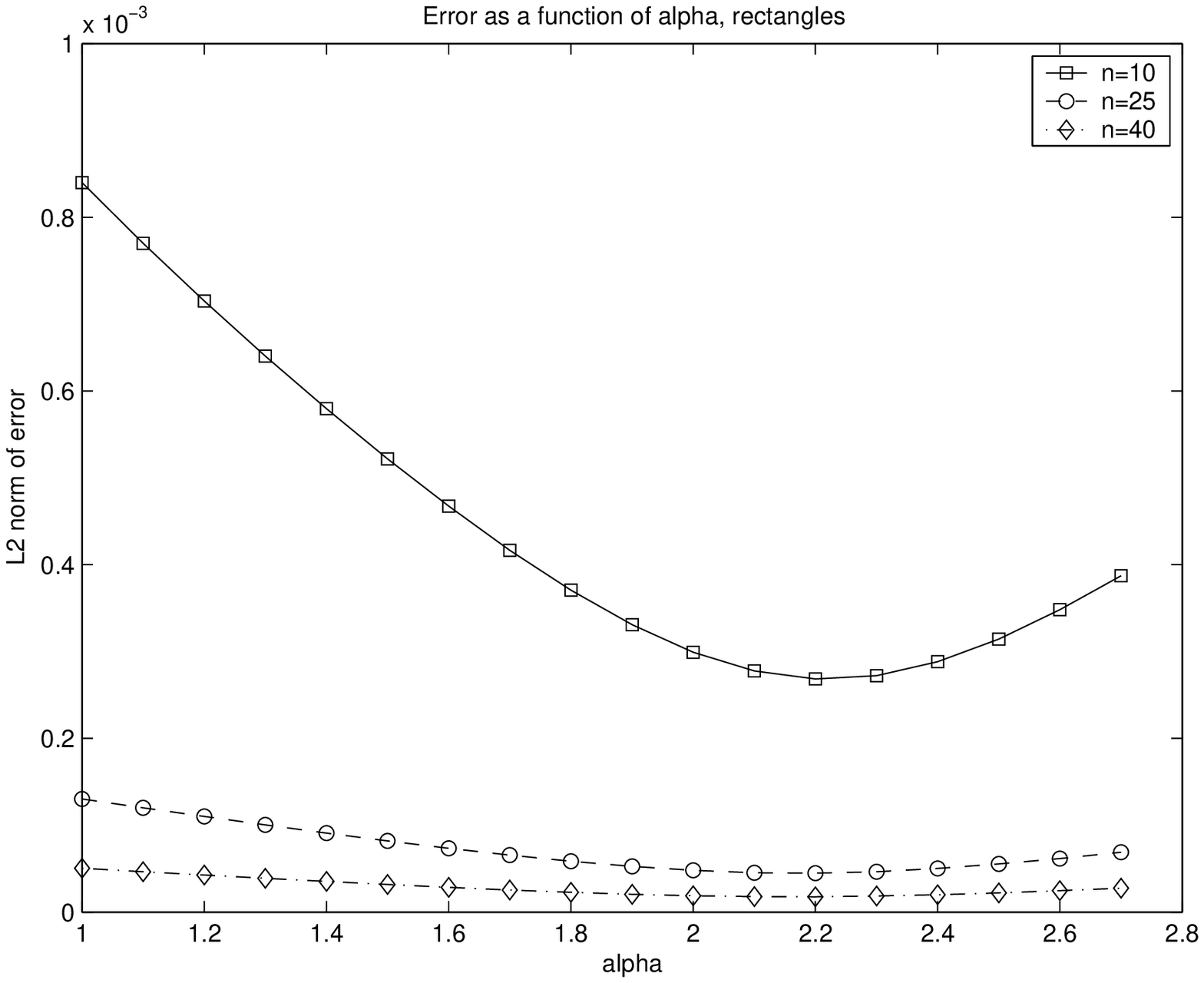}&
\includegraphics[width=0.45 \linewidth]{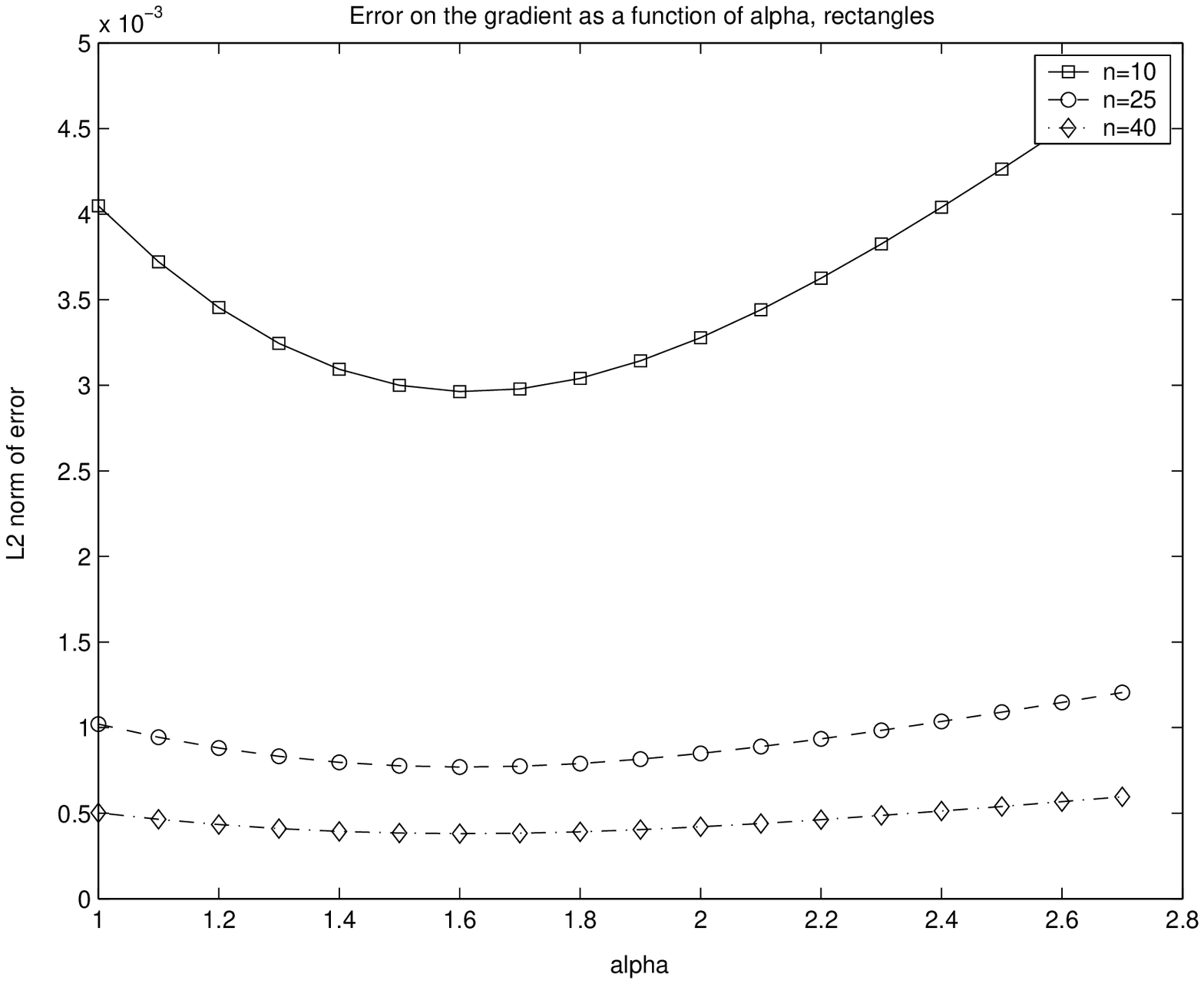}\end{tabular}
\end{center}\caption{Diagrams of the errors on the solution (left) and its 
gradient (right) for various sizes of 
 triangular (up) and rectangular (bottom) meshes, with respect to the value of the parameter $\alpha$}
 \label{alpha}
\end{figure}

Next, we tested different values of $\alpha$ to see how it affected the 
discretization error,  on the first anisotropic case. Although the value of 
$\alpha$ does influence the resulting discretization error, the optimal 
value seems to be independent on the mesh, in both the triangular and 
rectangular cases, see Figure \ref{alpha}. Note that in the case of the error on 
the solution itself, the numerical  optimal values for $\alpha$ are beyond the interval of 
convergence assumed in the theoretical analysis  $(0,1)$. 

These numerical tests therefore indicate that this use of 
an discrete gradient in finite volume schemes leads to
a correct numerical behavior, indeed comparable with low degree finite element schemes
on similar problems.

Finally, we replaced the point $x_K$ by the center of gravity of cell $K$
in the definition \refe{defaint},\refe{defaext} of the coefficients   $A_{K,L}$. In this case, we recall 
(see Remark \ref{rt}) that we obtain the discrete gradient based on the 
generalized Raviart-Thomas basis functions of \cite{convgrad}.
Indeed, the tests performed 
with this scheme for Case 1 or Case 2 did not yield correct approximations of 
the solution nor of its gradient.

\section{Conclusion}\label{conclu}
In this paper, we constructed a discrete gradient for piecewise constant functions. This discrete gradient   revealed several  advantages: it
is easy and cheap to compute, 
and it provides simple schemes for the approximation
of anisotropic 
diffusion convection problems. 
We showed a weak property convergence of this discrete gradient to the gradient of the
limit of the considered functions, together 
 with  a consistency property, both leading to the strong convergence
of the discrete solution and of its discrete gradient in the case
of a Dirichlet problem with full matrix diffusion.

Since this notion of admissible mesh includes  Vorono\"{\i} meshes,
which are more and more used in practice, and which seem to remain tractable even in
high space dimension, applications to financial mathematics problems are being studied \cite{berton}.
Applications to finite volume schemes for compressible Navier-Stokes equations are also expected to be 
succesful \cite{touazi}.
Further work includes a parametric study, and the generalization to
meshes without the orthogonality condition. 

\vspace{.5cm}
 
{\bf Aknowledgment} This work was supported by GDR MOMAS.


\begin{thebibliography}{00}

\bibitem{aav1}   I. Aavatsmark, T. Barkve, O. Boe and T. Mannseth,
Discretization on
unstructured grids for inhomogeneous, anisotropic media. Part I:
Derivation of the methods.  {\it SIAM Journal on Sc.  Comp.}, 19
(1998), 1700--1716.


\bibitem{aav2} I. Aavatsmark, T. Barkve, O. Boe and T. Mannseth,
Discretization on
unstructured grids for inhomogeneous, anisotropic media. Part II:
Discussion and numerical results.  {\it SIAM Journal on Sc.  Comp.}, 19 (1998), 1717--1736.

\bibitem{angermann} L. Angermann, 
A finite element method for the numerical solution of  
convection-dominated anisotropic diffusion equations.
{\it Numer. Math.}  85  (2000),  175--195.

\bibitem{berton} J. Berton
Comparaison de diff\'erentes m\'ethodes pour appr\'ecier les options am\'ericaines.
{\it Thesis of  Marne-la-Vall\'ee university (France)}, in preparation (2005).

 
\bibitem{letem} E. Ch\'enier, R. Eymard and X.  Nicolas,  
A Finite Volume Scheme for the Transport of Radionucleides 
{\it Porous Media: Simulation of Transport Around a Nuclear Waste Disposal Site: 
The COUPLEX Test Cases, Alain Bourgeat and Michel Kern eds,
Computational Geosciences}, 8 (2004), 163--172. 


\bibitem{coudiere}
Y. Coudi\`ere, J.P. Vila and P. Villedieu, 
Convergence rate of a finite volume scheme for a two-dimensional  
convection-diffusion problem.
{\it M2AN Math. Model. Numer. Anal.},  33  (1999), 493--516.

 
\bibitem{omnes} K. Domelevo, P. Omnes, 
 A finite volume method for the Laplace equation on almost arbitrary
 two-dimensional grids. {\it submitted} (2005). 




\bibitem{book}  R. Eymard, T. Gallou\"et and   R. Herbin,
Finite Volume Methods. {\it Handbook of Numerical Ana\-lysis, 
P.G. Ciarlet and J.L. Lions eds, North Holland}, 
7 (2000),  713--1020.

\bibitem{cvnl}  R. Eymard, T. Gallou\"et and   R. Herbin,
Convergence  of finite volume approximations to the solutions of
semilinear convection diffusion reaction equations. {\it Numer. Math.}, 82 (1999), 91--116.
 
\bibitem{convgrad}  R. Eymard, T. Gallou\"et and   R. Herbin,
Finite volume approximation of elliptic problems and convergence 
of an approximate gradient.
 {\it Appl. Num. Math.}, 37 (2001),31--53.

 
\bibitem{crasvf100}  R. Eymard, T. Gallou\"et and   R. Herbin,
A finite volume for anisotropic diffusion problems.
{\it  Comptes rendus \`a l'Acad\'emie des Sciences }, 339 (2004), 299--302.

\bibitem{convpardeg} R. Eymard, T. Gallou\"et, R. Herbin, A. Michel,
  Convergence of a finite volume scheme for nonlinear degenerate parabolic equations.
{\it  Num. Math},  92 (2002), 41--82.

\bibitem{voralhik} R. Eymard, D. Hilhorst,  M. Vohral\'{\i}k, 
Combined finite volume-nonconforming/mixed-hybrid finite element scheme 
for degenerate parabolic problems,
{\it  submitted} (2004).


\bibitem{f} { P.A. Forsyth,}   Control volume finite elements,  A control
volume finite element approach to NAPL
groundwater contamination. 
{\it  SIAM J. Sci. Stat. Comput.}, 12 (1991), 1029--1057.


  
\bibitem{fung}  L.S.-K. Fung,  L. Buchanan, and R. Sharma, 
Hybrid-CVFE Method for Flexible-
Grid Reservoir Simulation.
{\it   Soc. Pet. Eng. J.}, 19 (1994), 188-199.

\bibitem{ghv}
T. Gallou\"{e}t,  {  R. Herbin}  and  { M.H. 
Vignal},  Error estimate for the approximate finite volume solutions of convection 
diffusion 
equations with general boundary conditions.  {\it   SIAM J. 
Numer. Anal.}, 37 (2000), 1935--1972.  



\bibitem{vf4}  {R. Herbin},
An error estimate for a finite volume scheme for a
diffusion-convection problem on a triangular mesh. {\it  Num. Meth.
P.D.E.}
{  11} (1995),
165-173.
 


\bibitem{fvca1}  { R.  Herbin}, Finite
volume methods for
diffusion convection
equations  on general meshes. {\it in Finite volumes for complex
applications, Problems and Perspectives, F. Benkhaldoun and R.
Vilsmeier eds,
Hermes}, (1996) 153--160.

\bibitem{hermeline} F. Hermeline,  A finite volume method for 
the approximation of diffusion operators on  distorted meshes.
{\it  J. Comput. Phys. } 160  (2000),  481--499

  
   
\bibitem{nic2} X.H. Hu and R.A.  Nicolaides,  Covolume techniques for anisotropic media.
 {\it  Numer. Math.}  61  (1992),  215--234.
 
\bibitem{pt}  P. A. Jayantha   and Ian W. Turner,  
  A Second Order Finite Volume Technique for Simulating 
  Transport in Anisotropic Media",  {\it The Int. J. 
  of Num. Met. for Heat and Fluid Flow}, 13 (2003), 31--56.  
  
\bibitem{jt}  P. A. Jayantha and I. W. Turner, 
   A Second Order Control-Volume Finite-Element Least-Squares Strategy 
  for Simulating Diffusion in Strongly Anisotropic Media. {\it J.Comp. Math.},
   23 (2005), 1--16.  

\bibitem{LL} D. Lamberton and B. Lapeyre, 
An Introduction to Stochastic Calculus Applied to Finance.
{\it Chapman and Hall}, (1995).


 
\bibitem{mishev}
I.D.  Mishev,   Finite volume methods on Vorono\"{\i} meshes.
{ \it Num. Meth. P.D.E.}, 14 (1998),  193--212.


 
 \bibitem{nic1} R.A. Nicolaides,   Direct discretization of planar div-curl problems.
{\it SIAM J. Numer. Anal. }, 29  (1992), 32--56.
 

\bibitem{touazi} O. Touazi,
Mise en oeuvre d'un sch\'ema de volumes finis pour les \'equations de Navier-Stokes compressibles.
{\it Thesis of  Marne-la-Vall\'ee university (France)}, in preparation (2007).
  
\bibitem{putti} M. Putti and C. Cordes 
Finite Element Approximation of the Diffusion Operator on Tetrahedra.
{\it SIAM Journal on Scientific Computing}
19 (1998), 1154--1168.



\bibitem{wang} S. Wang, 
  Solving convection-dominated anisotropic diffusion equations by an  
exponentially fitted finite volume method.
{\it  Comput. Math. Appl.}  44  (2002), 1249--1265.
\end{thebibliography}
\end{document}